\documentclass{amsart}
\usepackage{amssymb}

\newtheorem{theorem}{Theorem}[section]
\theoremstyle{plain}
\newtheorem{lemma}[theorem]{Lemma}
\newtheorem{corollary}[theorem]{Corollary}
\newtheorem{proposition}[theorem]{Proposition}
\theoremstyle{definition}
\newtheorem{definition}[theorem]{Definition}
\newtheorem{example}[theorem]{Example}

\theoremstyle{remark}
\newtheorem{remark}[theorem]{Remark}
\numberwithin{equation}{section}

\input{tcilatex}

\begin{document}
\title[Loops and Semidirect Products]{Loops and Semidirect Products}
\author{Michael K. Kinyon}
\address{Department of Mathematics \& Computer Science\\
Indiana University South Bend\\
South Bend, IN 46634 USA}
\email{mkinyon@iusb.edu}
\urladdr{http://www.iusb.edu/\symbol{126}mkinyon}
\author{Oliver Jones}
\address{Department of Mathematics\\
University of Puerto Rico\\
Mayaguez Campus\\
PO Box 9018\\
Mayaguez, PR 00681-9018}
\email{o\_jones@rumac.upr.clu.edu}
\maketitle


\section{Introduction}

A \emph{left loop} $(B,\cdot )$ is a set $B$ together with a binary
operation $\cdot $ such that (i) for each $a\in B$, the mapping $%
x\longmapsto a\cdot x$ is a bijection, and (ii) there exists a two-sided
identity $1\in B$ satisfying $1\cdot x=x\cdot 1=x$ for every $x\in B$. A
right loop is similarly defined, and a \emph{loop} is both a right loop and
a left loop \cite{bruck} \cite{cps}.

In this paper we study semidirect products of loops with groups. This is a
generalization of the familiar semidirect product of groups. Recall that if $%
G$ is a group with subgroups $B$ and $H$ where $B$ is normal, $G=BH$, and $%
B\cap H=\left\{ 1\right\} $, then $G$ is said to be an \emph{internal}
semidirect product of $B$ with $H$. On the other hand, if $B$ and $H$ are
groups and $\sigma :H\rightarrow \mathrm{Aut}(B):h\mapsto \sigma _{h}$ is a
homomorphism, then the \emph{external} semidirect product of $B$ with $H$
given by $\sigma $, denoted $B\rtimes _{\sigma }H$, is the set $B\times H$
with the multiplication 
\begin{equation}
(a,h)(b,k)=(a\cdot \sigma _{h}(b),hk).  \label{extern-semi}
\end{equation}
A special case of this is the \emph{standard} semidirect product where $H$
is a subgroup of the automorphism group of $B$, and $\sigma $ is the
inclusion mapping. The relationship between internal, external and standard
semidirect products is well known.

These considerations can be generalized to loops. We now describe the
contents of the sequel.

In \S 2, we consider the natural embedding of a left loop $B$ into its
permutation group $\mathrm{Sym}(B)$. This leads to a factorization of $%
\mathrm{Sym}(B)$ into a sub\emph{set} $L(B)$ consisting of the left
translations of $B$ and a sub\emph{group} $\mathrm{Sym}_{1}(B)$ consisting
of permutations fixing the identity element $1\in B$. We then discuss left
inner mappings and deviations \cite{sabinin}, and show how these
characterize those permutations which are pseudo-automorphisms and
automorphisms. We discuss how the aforementioned factorization of $\mathrm{%
Sym}(B)$ is related to the group multiplication; here the left inner
mappings and deviations play a role in decomposing the product of
permutations. Finally, we give Sabinin's definition of the \emph{standard}
semidirect product of a left loop $B$ with one of its transassociants \cite
{sabinin}. This semidirect product has occasionally been rediscovered for
various classes of loops. We conclude the section with an example.

In \S 3, we consider \emph{internal} semidirect products of left loops and
groups: given a group $G$, a subgroup $H<G$, and a transversal $B\subseteq G$
of $H$ which contains the identity, $B$ naturally has the structure of a
left loop. This is equivalent to putting a loop structure on the set $G/H$
of cosets \cite{baer}, but it is closer to the examples to work with
transversals. We consider the relationship between the loop structure of $B$
and the multiplication in $G$, paralleling the discussion in \S 2. We follow
Sabinin \cite{sabinin}, but with considerably more detail; see also \cite
{kiechle}. We then consider conditions under which subgroups and factor
groups inherit the semidirect product structure. A particular case of the
latter is obtained by modding out the core of $H$ in $G$.

In \S 4, we use transversal decompositions $G=BH$ to study certain loop
identities, especially those related to Bol loops and Bruck loops. This part
of our study is related to work of Ungar \cite{u-wags} and Kreuzer and
Wefelscheid \cite{kreuz-wefel}, but we do not assume as much structure at
the outset. We then give examples of internal semidirect products,
illustrating some of the results.

In \S 5, we generalize the standard semidirect product of a left loop $B$
with a particular subgroup of $\mathrm{Sym}_{1}(B)$ to an \emph{external}
semidirect product of a left loop $B$ with a group $H$. As for the usual
semidirect product of groups, the main interest here is in the case where
the defining homomorphism from $H$ to $\mathrm{Sym}_{1}(B)$ is not
injective. Our construction seems to be new, and we give examples. We
conclude by discussing how the three semidirect products are related,
generalizing the relationship between the usual standard, internal, and
external semidirect products of groups.

There exist notions of semidirect products of loops which are different from
that which we consider here. One definition is as follows: a loop $R$ is an
internal semidirect product of the normal subloop $P$ by the subloop $Q$ if $%
R=PQ$ and $P\cap Q=\left\{ 1\right\} $. This definition was given by
Birkenmeier \emph{et al} \cite{birk1} and Birkenmeier and Xiao \cite{birk2},
who studied nonassociative loops which are internal semidirect products of
groups. Goodaire and Robinson \cite{good-robin} defined an internal
semidirect product similarly with additional conditions given in terms of
associators. In contrast, our internal semidirect product follows Sabinin 
\cite{sabinin}: $G=BH$ is a factorization of a group $G$ into a subgroup $H$
and a transversal $B$. Even if $B$ with its induced operation turns out to
be a group, it is not necessarily a subgroup of $G$. Also, $H$ does not
necessarily stabilize $B$ by conjugation. Thus the two notions of semidirect
product are quite distinct.

In group theory, the question of which groups have a semidirect product
structure is answered by cohomology theory. Cohomology has been generalized
to loops in at least two distinct ways; see Eilenberg and MacLane \cite
{maclane}, and Johnson and Leedham-Green \cite{johnson}. At present, we do
not know if the semidirect product of the present paper has a suitable
cohomological interpretation.

\section{Standard Semidirect Products}

For a set $B$, let $\mathrm{Sym}(B)$ denote the group of permutations of $B$%
. If $(B,\cdot )$ is a left loop with identity element $1\in B$, let $%
\mathrm{Sym}_{1}(B)$ denote the subgroup consisting of permutations fixing $%
1 $. For each $x\in B$, define the left translation mapping $%
L_{x}:B\rightarrow B$ by $L_{x}y=x\cdot y$. Define the left division
operation $\backslash :B\times B\rightarrow B$ by $x\backslash y=L_{x}^{-1}y$
for $x,y\in B$. We denote the right inverse of $x\in B$ by $x^{\rho
}=x\backslash 1$. The right inversion mapping $\rho :x\longmapsto x^{\rho }$
fixes $1$, and is a permutation if and only if each $x\in B$ has a unique
left inverse $x^{\lambda }$.

Let $L(B)=\left\{ L_{x}:x\in B\right\} $ be the set of left translations,
and let $\mathrm{LMlt}(B,\cdot )=\langle L(B)\rangle $ be the \emph{left
multiplication group}, i.e., $\mathrm{LMlt}(B,\cdot )$ is the subgroup of $%
\mathrm{Sym}(B)$ generated by $L(B)$. Let $\mathrm{LMlt}_{1}(B,\cdot )=%
\mathrm{LMlt}(B,\cdot )\cap \mathrm{Sym}_{1}(B)$. For $x,y\in B$, the
permutation $L(x,y)\in \mathrm{Sym}(B)$ defined by 
\begin{equation}
L(x,y)=L_{xy}^{-1}L_{x}L_{y}  \label{eq:left_inner}
\end{equation}
is called a \emph{left inner mapping}. Clearly 
\begin{equation}
L(1,x)=L(x,1)=I  \label{eq:L(1,x)=L(x,1)=I}
\end{equation}
for all $x\in B$. We have $L(x,y)\in \mathrm{LMlt}_{1}(B,\cdot )$, and in
fact, 
\begin{equation}
\mathrm{LMlt}_{1}(B,\cdot )=\langle L(x,y):x,y\in B\rangle
\label{eq:LMlt_generators}
\end{equation}
(\cite{bruck}, p.61; \cite{pflugfelder}, I.5.2). $\mathrm{LMlt}_{1}(B,\cdot
) $ is called the \emph{left inner mapping group}.

A permutation $\phi \in \mathrm{Sym}(B)$ is called a \emph{%
pseudo-automorphism} with \emph{companion} $c\in B$ if $c\cdot \phi (x\cdot
y)=(c\cdot \phi (x))\cdot \phi (y)$ for all $x,y\in B$, or equivalently, 
\begin{equation}
L_{c}\phi L_{x}\phi ^{-1}=L_{c\cdot \phi (x)}  \label{eq:pseudo-aut}
\end{equation}
for all $x\in B$. The set $\mathrm{PsAut}(B)$ of pseudo-automorphisms of $B$
is a group under composition of mappings. Since a pseudo-automorphism fixes $%
1$ (take $y=1$ and cancel $c\cdot \phi (x)$), $\mathrm{PsAut}(B)$ is a
subgroup of $\mathrm{Sym}_{1}(B)$. A left loop $(B,\cdot )$ is said to have
the \emph{pseudo-}$A_{l}$ \emph{property} if every left inner mapping $%
L(x,y) $ is a pseudo-automorphism. By (\ref{eq:LMlt_generators}), this is
equivalent to the assertion that $\mathrm{LMlt}_{1}(B,\cdot )\leq \mathrm{%
PsAut}(B)$.

A pseudo-automorphism with companion $1$ is an \emph{automorphism} of $B$.
Let $\mathrm{Aut}(B)$ denote the group of automorphisms of $B$. A left loop $%
(B,\cdot )$ is said to have the $A_{l}$ \emph{property} if every left inner
mapping $L(x,y)$ is an automorphism. By (\ref{eq:LMlt_generators}), this is
equivalent to the assertion that $\mathrm{LMlt}_{1}(B,\cdot )\leq \mathrm{Aut%
}(B)$.

For each $x\in B$ and $\phi \in \mathrm{Sym}(B)$, the permutation $\mu
_{x}(\phi )\in \mathrm{Sym}(B)$ defined by 
\begin{equation}
\mu _{x}(\phi )=L_{\phi (x)}^{-1}\phi L_{x}\phi ^{-1}  \label{eq:deviation}
\end{equation}
is called the \emph{deviation} of $\phi $ at $x$ (\cite{sabinin} \cite
{kiechle}, 2.C). Clearly, we have 
\begin{eqnarray}
\mu _{1}(\phi ) &=&I  \label{eq:mu(1,phi)=I} \\
\mu _{x}(I) &=&I  \label{eq:mu(x,I)=I}
\end{eqnarray}
for all $x\in B$, $\phi \in \mathrm{Sym}(B)$. As the next result shows,
deviations measure how much arbitrary permutations ``deviate'' from being
(pseudo-)automorphisms.

\begin{proposition}
\label{prop:deviation_properties}Let $(B,\cdot )$ be a left loop, and let $%
\phi \in \mathrm{Sym}(B)$ be given.

\begin{enumerate}
\item  $\phi \in \mathrm{Sym}_{1}(B)$ if and only if $\mu _{x}(\phi )\in 
\mathrm{Sym}_{1}(B)$ for all $x\in B$.

\item  $\phi \in \mathrm{PsAut}(B)$ if and only if there exists $c\in B$
such that $\mu _{x}(\phi )=L(c,\phi (x))^{-1}$ for all $x\in B$.

\item  $\phi \in \mathrm{Aut}(B)$ if and only if $\mu _{x}(\phi )=I$ for all 
$x\in B$.
\end{enumerate}
\end{proposition}

\begin{proof}
1. For $x\in B$, $\mu _{x}(\phi )1=1$ if and only if $\phi (x)=\phi (x\cdot
\phi ^{-1}(1))$ if and only if $x=x\cdot \phi ^{-1}(1)$ if and only if $%
1=\phi ^{-1}(1)$.

\noindent 2. For all $x,c\in B$, a computation using (\ref{eq:deviation})
and (\ref{eq:left_inner}) gives 
\begin{equation}
L_{c}\phi L_{x}\phi ^{-1}=L_{c\cdot \phi (x)}L(c,\phi (x))\mu _{x}(\phi ).
\label{eq:dev_tmp}
\end{equation}
From this, it is clear that $\phi $ is a pseudo-automorphism with companion $%
c$ if and only if $L(c,\phi (x))\mu _{x}(\phi )=I$.

\noindent 3. This follows from taking $c=1$ in (\ref{eq:dev_tmp}) and using (%
\ref{eq:L(1,x)=L(x,1)=I}).
\end{proof}

For a left loop $(B,\cdot )$, the mapping $B\rightarrow L(B):x\longmapsto
L_{x}$ is bijective ($L_{x}=L_{y}$ implies $x=L_{x}(1)=L_{y}(1)=y$). Note
that $L(B)$ itself can be given the structure of a left loop isomorphic to $%
(B,\cdot )$ with the obvious definition: 
\begin{equation}
L_{x}\cdot L_{y}=L_{x\cdot y}  \label{eq:L(B)_op}
\end{equation}
for $x,y\in B$.

Let $G$ be any group satisfying $\mathrm{LMlt}(B,\cdot )\leq G\leq \mathrm{%
Sym}(B)$, and let $H=G\cap \mathrm{Sym}_{1}(B)$, so that $\mathrm{LMlt}%
_{1}(B,\cdot )\leq H\leq \mathrm{Sym}_{1}(B)$. For any $\phi \in G$, we have 
$\phi =L_{x}\psi $ where $x=\phi (1)$ and $\psi =L_{x}^{-1}\phi $. Clearly $%
\psi (1)=1$, and thus $\psi \in H$ since $G$ contains $\mathrm{LMlt}(B,\cdot
)$. The factorization of $\phi $ into a left translation $L_{x}$ in $L(B)$
and a permutation $\psi $ in $H$ is unique. Indeed, if $L_{x}\psi
=L_{y}\varphi $ for $x,y\in B$, $\psi ,\varphi \in H$, then applying both
sides to $1$ gives $x=y$, and thus $L_{x}\psi =L_{x}\varphi $; cancelling $%
L_{x}$ gives $\psi =\varphi $.

Summarizing, for any group $G$ satisfying $\mathrm{LMlt}(B,\cdot )\leq G\leq 
\mathrm{Sym}(B)$, we have the following decomposition: 
\begin{equation}
G=L(B)H  \label{eq:G=L(B)H}
\end{equation}
where $H=G\cap \mathrm{Sym}_{1}(B)$. The factorization of elements is
unique, and we also have 
\begin{equation}
L(B)\cap H=\{I\}.  \label{eq:L(B)capH=I}
\end{equation}

It is natural to ask how the factorization (\ref{eq:G=L(B)H}) of a group $G$
into a sub\emph{set} $L(B)$ with a left loop structure given by (\ref
{eq:L(B)_op}) and a sub\emph{group} $H$ interacts with the multiplication in 
$G$. We first examine this question for $G=\mathrm{Sym}(B)$ and $H=\mathrm{%
Sym}_{1}(B)$. For permutations $L_{x}\phi $ and $L_{y}\psi $ with $x,y\in B$%
, $\phi ,\psi \in \mathrm{Sym}_{1}(B)$, we have $(L_{x}\phi L_{y}\psi
)(1)=x\cdot \phi (y)$. Also, $L_{x\cdot \phi (y)}^{-1}L_{x}\phi
L_{y}=L(x,\phi (y))\mu _{y}(\phi )\phi $. Put together, these observations
give the following factorization of a product in $\mathrm{Sym}(B)=L(B)\ 
\mathrm{Sym}_{1}(B)$.

\begin{proposition}
\label{prop:factored_product}Let $(B,\cdot )$ be a left loop. For all $%
x,y\in B$, $\phi ,\psi \in \mathrm{Sym}_{1}(B)$, 
\begin{equation}
(L_{x}\phi )(L_{y}\psi )=L_{x\cdot \phi (y)}\left[ L(x,\phi (y))\mu
_{y}(\phi )\phi \psi \right] .  \label{fac-prod}
\end{equation}
\end{proposition}

Using (\ref{eq:L(B)_op}), (\ref{fac-prod}) can be rewritten as 
\begin{equation*}
(L_{x}\phi )(L_{y}\psi )=(L_{x}\cdot L_{\phi (y)})(L(x,\phi (y))\mu
_{y}(\phi )\phi \psi ).
\end{equation*}
Thus we see that the operation $L_{x}\cdot L_{y}$ defined by (\ref
{eq:L(B)_op}) is simply the projection of the composition $L_{x}L_{y}$ onto $%
L(B)$.

For the factorization (\ref{fac-prod}) to hold in a subgroup $G$ of $\mathrm{%
Sym}(B)$, it is clear that it is necessary that the part of (\ref{fac-prod})
in square brackets be in the subgroup $H=G\cap \mathrm{Sym}_{1}(B)$. Thus
assume $L(x,\phi (y))\mu _{y}(\phi )\in H$ for all $x,y\in B$, $\phi \in H$.
Taking $\phi =I$, we have $L(x,y)\in H$ for all $x,y\in H$, and this implies 
$\mu _{y}(\phi )\in H$ for all $y\in B$, $\phi \in H$. This leads us to the
following definition of Sabinin \cite{sabinin}. A subgroup $H\leq \mathrm{Sym%
}_{1}(B)$ is said to be a \emph{transassociant} of $B$ if $L(x,y)\in H$ for
all $x,y\in B$ and if $\mu _{x}(\phi )\in H$ for all $x\in B$, $\phi \in H$.

\begin{proposition}
\label{prop:transassociants}Let $(B,\cdot )$ be a left loop.

\begin{enumerate}
\item  $\mathrm{LMlt}_{1}(B,\cdot )$ is a transassociant.

\item  If $(B,\cdot )$ has the pseudo-$A_{l}$ property, and if $\mathrm{LMlt}%
_{1}(B,\cdot )\leq H\leq \mathrm{PAut}(B)$, then $H$ is a transassociant.

\item  If $(B,\cdot )$ has the $A_{l}$ property, and if $\mathrm{LMlt}%
_{1}(B,\cdot )\leq H\leq \mathrm{Aut}(B)$, then $H$ is a transassociant.
\end{enumerate}
\end{proposition}

\begin{proof}
\noindent 1. If $\phi \in \mathrm{LMlt}_{1}(B,\cdot )$, then by (\ref
{eq:deviation}) and Proposition \ref{prop:deviation_properties}(1), $\mu
_{x}(\phi )\in \mathrm{LMlt}_{1}(B,\cdot )$ for all $x\in B$.

\noindent 2. This follows from Proposition \ref{prop:deviation_properties}%
(2).

\noindent 3. This follows from Proposition \ref{prop:deviation_properties}%
(3).
\end{proof}

Consider again the case where $G=\mathrm{Sym}(B)=L(B)\ \mathrm{Sym}_{1}(B)$.
The factorization (\ref{eq:G=L(B)H}), along with (\ref{eq:L(B)capH=I}),
gives a one-to-one correspondence between $B\times \mathrm{Sym}_{1}(B)$ and $%
\mathrm{Sym}(B)$ given by $(x,\phi )\longmapsto L_{x}\phi $. Thus we may use
Proposition \ref{prop:factored_product} to define a binary operation on $%
B\times \mathrm{Sym}_{1}(B)$: 
\begin{equation}
(x,\phi )\cdot (y,\psi )=(x\cdot \phi (y),L(x,\phi (y))\mu _{y}(\phi )\phi
\psi )  \label{semi-op}
\end{equation}
for $x,y\in B$, $\phi ,\psi \in \mathrm{Sym}_{1}(B)$. By construction, $%
(B\times \mathrm{Sym}_{1}(B),\cdot )$ is a group isomorphic to $\mathrm{Sym}%
(B)$. We now present Sabinin's definition \cite{sabinin} of the semidirect
product of a left loop with one of its transassociants.

\begin{definition}
\label{defn-semidir-product}Let $B$ be a left loop and let $H\leq \mathrm{Sym%
}_{1}(B)$ be a transassociant of $B$. Define a binary operation $\cdot $ on
the set $B\times H$ as follows: 
\begin{equation*}
(x,\phi )\cdot (y,\psi )=(x\cdot \phi (y),L(x,\phi (y))\mu _{y}(\phi )\phi
\psi )
\end{equation*}
for all $a,b\in B$, $\phi ,\psi \in H$. Then $(B\times H,\cdot )$ is called
the \emph{standard} \emph{semidirect product} of $B$ with $H$, and is
denoted $B\rtimes H$.
\end{definition}

It is immediate from this definition that $B\rtimes H$ is a subgroup of the
group $B\rtimes \mathrm{Sym}_{1}(B)\cong \mathrm{Sym}(B)$. We have proven
the following result.

\begin{proposition}
(\cite{sabinin}, Thm. 2) Let $B$ be a left loop and let $H\leq \mathrm{Sym}%
_{1}(B)$ be a transassociant of $B$. Then $B\rtimes H$ is a group.
\end{proposition}

Incidentally, as Sabinin \cite{sabinin} has noted, in order for $B\rtimes H$
to be a group, it is only necessary for $B$ to have a \emph{right} identity
element $1$. In this case, (\ref{eq:L(B)capH=I}) becomes $L(B)\cap
H=\{L_{1}\}$ because $L_{1}\neq I$.

We now consider some special cases.

\begin{remark}
\label{rem:special_standard}

\begin{enumerate}
\item  If $(B\cdot )$ is a group, then the product in $B\rtimes H$ is given
by 
\begin{equation*}
(x,\phi )\cdot (y,\psi )=(x\cdot \phi (y),\mu _{y}(\phi )\phi \psi ),
\end{equation*}
$x,y\in B$, $\phi ,\psi \in H$. This generalized semidirect product of
groups was rediscovered by Jajcay \cite{jajcay}, who dubbed it the
``rotary'' product of groups. The semidirect product $\mathrm{Sym}(B)\cong
B\rtimes \mathrm{Sym}_{1}(B)$ can be seen as a detailed description of the
algebraic structure of the regular representation of $(B\cdot )$.

\item  Assume $(B\cdot )$ is a pseudo-$A_{l}$ left loop and that $\mathrm{%
LMlt}_{1}(B,\cdot )\leq H\leq \mathrm{PsAut}(B)$. By Proposition \ref
{prop:transassociants}(2), $H$ is a transassociant. In this case, the
product in $B\rtimes H$ is given by 
\begin{equation*}
(x,\phi )\cdot (y,\psi )=(x\cdot \phi (y),L(x,\phi (y))L(c,\phi
(y))^{-1}\phi \psi ),
\end{equation*}
$x,y\in B$, $\phi ,\psi \in H$, where $c$ is a companion of $\phi $, using
Proposition \ref{prop:deviation_properties}(2). The semidirect product group 
$\mathrm{PsAff}(B)\cong B\rtimes \mathrm{PsAut}(B)$ is called the \emph{%
pseudo-affine group} of $(B\cdot )$.

\item  Assume $(B\cdot )$ is an $A_{l}$ left loop and that $\mathrm{LMlt}%
_{1}(B,\cdot )\leq H\leq \mathrm{Aut}(B)$. By Proposition \ref
{prop:transassociants}(3), $H$ is a transassociant. In this case, the
product in $B\rtimes H$ is given by 
\begin{equation*}
(x,\phi )\cdot (y,\psi )=(x\cdot \phi (y),L(x,\phi (y))\phi \psi ),
\end{equation*}
$x,y\in B$, $\phi ,\psi \in H$, using Proposition \ref
{prop:deviation_properties}(3). The semidirect product group $\mathrm{Aff}%
(B)\cong B\rtimes \mathrm{Aut}(B)$ is called the \emph{affine group} of $%
(B,\cdot )$. For $A_{l}$ left loops with the left inverse property (see \S
4), this semidirect product was rediscovered by Kikkawa \cite{kikkawa} and
later, using different terminology, by Ungar \cite{u-wags}.

\item  If $B$ is a group and $H$ is a subgroup of $\mathrm{Aut}(B)$, then $%
B\rtimes H$ is the usual standard semidirect product of groups.
\end{enumerate}
\end{remark}

We now give an explicit example of a standard semidirect product.

\begin{example}
\label{ex-disk}Let $\mathbb{D}=\{z\in \mathbb{C}:\left| z\right| <1\}$ be
the complex unit disk. Note that the circle group $S^{1}=\left\{ a\in 
\mathbb{C}:\left| a\right| =1\right\} $ acts on $\mathbb{D}$ by
multiplication of complex numbers. For $x,y\in \mathbb{D}$, define 
\begin{equation}
x\oplus y=\frac{x+y}{1+\bar{x}y}\text{.}  \label{disk-op}
\end{equation}
$(\mathbb{D},\oplus )$ turns out to be a $B$-loop \cite{k-u} \cite{u-disk}
(see \S 3 for the definition). The left inner mappings are given by
unimodular complex numbers 
\begin{equation}
L(x,y)z=\frac{1+x\bar{y}}{1+\bar{x}y}z  \label{disk-gyr}
\end{equation}
for $x,y,z\in \mathbb{D}$. If we identify $S^{1}$ with its natural image in $%
\mathrm{Aut}(\mathbb{D})$, then $\mathrm{Aut}(\mathbb{D})$ is generated by $%
S^{1}$ and the complex conjugation mapping $x\longmapsto \bar{x}$. The
semidirect product $\mathbb{D}\rtimes S^{1}$ turns out to be isomorphic to
the orientation-preserving M\"{o}bius group of the disk $\mathbb{D}$. The
semidirect product $\mathbb{D}\rtimes \mathrm{Aut}(\mathbb{D})$ is
isomorphic to the full M\"{o}bius group of $\mathbb{D}$.
\end{example}

\section{Internal Semidirect Products}

Let $G$ be a group with identity element $1$, let $H$ be a subgroup of $G$,
and let $B$ be a left \emph{transversal} of $H$ in $G$, i.e., for every $%
g\in G$, there exists a unique $a\in B$ and a unique $h\in H$ such that $%
g=ah $. (Equivalently, each element of $B$ is a representative of a unique
left coset of $H$ in $G$.) We will call the factorization $G=BH$ a \emph{%
transversal decomposition}. (Sabinin \cite{sabinin} calls $B$ a
``quasi-reductant''. Here we adapt more standard group-theoretic
terminology.)

Let $e\in B$ denote the representative of the coset $1H=H$. Then obviously $%
B\cap H=\{e\}$ (since $eH=H$). We have $e^{-1}\in H$. Now let $\tilde{B}%
=Be^{-1}=\{ae^{-1}:a\in B\}$. Then $\tilde{B}$ is a transversal of $H$ in $G$%
. Indeed, let $g\in G$ satisfy $g=ah$ where $a\in B$, $h\in H$. Then $%
g=(ae^{-1})(eh)$. On the other hand, if $g=bk$ for $b\in \tilde{B}$, $k\in H$%
, then $be\in B$ and $ae^{-1}=g(eh)^{-1}=(be)(e^{-1}kh^{-1}e^{-1})$. Thus $%
a=be$ and $k=eh$, which shows uniqueness of the decomposition. In addition,
we clearly have $\tilde{B}\cap H=\{1\}$. This discussion shows that there is
no real loss in assuming that $e=1$, i.e., 
\begin{equation}
B\cap H=\{1\}.  \label{eq:BcapH=e}
\end{equation}
In this case, $B$ is called a \emph{unital} transversal and $G=BH$ is called
a \emph{unital} transversal decomposition. ($B$ is sometimes called a
``uniform quasi-reductant'' \cite{sabinin}, and $G=BH$ is called an
``exact'' decomposition \cite{kram-urb} \cite{kreuz-wefel} \cite{u-wags}.)
We will assume throughout that our transversal decompositions are unital
without specifically mentioning it.

Let $G=BH$ be a transversal decomposition. Define a binary operation $\cdot
:B\times B\rightarrow B$ as follows: for $x,y\in B$, $x\cdot y\in B$ is
defined by 
\begin{equation}
(x\cdot y)H=xyH.  \label{eq:B_operation}
\end{equation}
Also define 
\begin{equation}
l(x,y)=(x\cdot y)^{-1}xy.  \label{eq:l(a,b)}
\end{equation}
We call $l:B\times B\rightarrow H$ the \emph{transversal mapping}. Note that 
\begin{eqnarray}
1\cdot x &=&x\cdot 1=x  \label{eq:1_identity} \\
l(1,x) &=&l(x,1)=1  \label{eq:l(1,a)=l(a,1)=1}
\end{eqnarray}
for all $x\in B$. Finally we define $\backslash :B\times B\rightarrow B$ as
follows: for $x,y\in B$, $x\backslash y\in B$ is defined by 
\begin{equation}
(x\backslash y)H=x^{-1}yH  \label{eq:left_division}
\end{equation}
for $x,y\in B$.

The following result is well-known (e.g., \cite{sabinin}, Thm. 7; \cite
{kreuz-wefel}, Thm. 3.2; \cite{kiechle}, Thm. 2.7.).

\begin{proposition}
\label{prop:left_loop}$(B,\cdot )$ is a left loop.
\end{proposition}

\begin{proof}
From (\ref{eq:1_identity}), $1$ is a two-sided identity. For $x,y\in B$, we
have 
\begin{equation*}
(x\cdot (x\backslash y))H=x(x\backslash y)H=xx^{-1}yH=yH
\end{equation*}
and 
\begin{equation*}
(x\backslash (x\cdot y))H=x^{-1}(x\cdot y)H=x^{-1}xyH=yH.
\end{equation*}
Thus each left translation $L_{x}:B\rightarrow B:y\longmapsto x\cdot y$ has
an inverse given by $L_{x}^{-1}(y)=x\backslash y$.
\end{proof}

Obviously $(B,\cdot )$ is a subgroup of $G$ if and only if the transversal
mapping $l:B\times B\rightarrow H$ is trivial, i.e., $l(x,y)=1$ for all $%
x,y\in B$.

For $x\in B$, recall that $x^{\rho }$ is the unique right inverse of $x$,
i.e., $x\cdot x^{\rho }=1$ or $x^{\rho }=x\backslash 1$. Thus $x^{\rho
}H=x^{-1}H$, i.e., $x^{\rho }$ is the representative of $x^{-1}H$ in $B$. In
addition, we have 
\begin{equation}
l(x,x^{\rho })=xx^{\rho }.  \label{eq:l(a,ar)=aar}
\end{equation}
In general, the mapping $\rho :B\rightarrow B:x\longmapsto x^{\rho }$ is not
a permutation, i.e., not every element has a unique left inverse. The next
result characterizes left loops for which this holds.

\begin{proposition}
\label{prop:left_inverse}Let $G=BH$ be a transversal decomposition. The
following are equivalent: (i) $B$ is a right transversal of $H$ in $G$; (ii) 
$B^{-1}$ is a left transversal of $H$ in $G$; (iii) each element of $B$ has
a unique left inverse in $(B,\cdot )$.
\end{proposition}

\begin{proof}
The equivalence of (i) and (ii) is obvious. Since $B$ is a transversal, $%
B^{-1}$ is a transversal if and only if, for every $x\in B$, there exists a
unique $x^{\lambda }\in B$ such that $xH=(x^{\lambda })^{-1}H$. This is
equivalent to $x^{\lambda }xH=H$, or $x^{\lambda }\cdot x=1$, which
establishes the equivalence of (ii) and (iii).
\end{proof}

There is a natural action of $G$ on $B$ which is given by the action of $G$
on the set of left cosets $G/H$. For $g\in G$, $x\in B$, we define $\sigma
_{g}(x)\in B$ by 
\begin{equation}
\sigma _{g}(x)H=gxH.  \label{eq:sigma}
\end{equation}
for all $x\in B$, $h\in H$. Recall that the \emph{core} of the subgroup $H$
is defined by 
\begin{equation}
\mathrm{core}_{G}(H)=\bigcap_{g\in G}gHg^{-1},  \label{eq:core}
\end{equation}
i.e., $\mathrm{core}_{G}(H)$ is largest subgroup of $H$ which is normal in $%
G $. For the next result, see also \cite{kiechle}, Thm. 2.8.

\begin{theorem}
\label{thm:sigma_properties}Let $G=BH$ be a transversal decomposition.

\begin{enumerate}
\item  $\sigma :G\rightarrow \mathrm{Sym}(B)$ is a homomorphism.

\item  $\sigma (H)\leq \mathrm{Sym}_{1}(B)$.

\item  For all $x,y\in B$, $h\in H$, 
\begin{eqnarray}
\sigma _{x} &=&L_{x}  \label{eq:sigma_a=L_a} \\
\sigma _{l(x,y)} &=&L(x,y)  \label{eq:sigma_l(a,b)=L(a,b)}
\end{eqnarray}

\item  $\ker (\sigma )=\mathrm{core}_{G}(H)$.
\end{enumerate}
\end{theorem}

\begin{proof}
1. For $g,g^{\prime }\in G$, $x\in B$, we have by (\ref{eq:sigma}), 
\begin{equation*}
(\sigma _{g}\sigma _{g^{\prime }})(x)H=g\sigma _{g^{\prime }}(x)H=gg^{\prime
}xH=\sigma _{gg^{\prime }}(x)H.
\end{equation*}
In particular, $\sigma _{g^{-1}}\sigma _{g}=\sigma _{g}\sigma _{g^{-1}}=I$,
and thus $\sigma (G)\leq \mathrm{Sym}(B)$.

\noindent 2. For $h\in H$, $x\in B$, we have by (\ref{eq:sigma}), 
\begin{equation*}
\sigma _{h}(1)H=h1H=H=1H,
\end{equation*}
and thus $\sigma _{h}(1)=1$ as claimed.

\noindent 3. For $x,y\in B$, (\ref{eq:sigma}) and (\ref{eq:B_operation})
imply $\sigma _{x}(y)=x\cdot y=L_{x}y$, which gives (\ref{eq:sigma_a=L_a}).
For $x,y,z\in B$, we use (\ref{eq:sigma}), (\ref{eq:l(a,b)}), (\ref
{eq:B_operation}), and (\ref{eq:left_inner}) to compute 
\begin{eqnarray*}
\sigma _{l(x,y)}(z)H &=&l(x,y)zH=(x\cdot y)^{-1}xyzH \\
&=&[(x\cdot y)^{-1}\cdot (x\cdot (y\cdot z))]H \\
&=&[(L_{x\cdot y}^{-1}L_{x}L_{y})(z)]H \\
&=&(L(x,y)z)H.
\end{eqnarray*}
This establishes (\ref{eq:sigma_l(a,b)=L(a,b)}).

\noindent 4. If $h\in \ker (\sigma )$, then for all $x\in B$, $\sigma
_{h}(x)H=hxH=xH$, i.e., for all $x\in B$, $hx\in xH$. This is equivalent to $%
hg\in gH$ for all $g\in G$, or $h\in \bigcap_{g\in G}gHg^{-1}$, as claimed.
\end{proof}

For $x\in B$, $h\in H$, define $m(x,h)\in H$ by 
\begin{equation}
m(x,h)=\sigma _{h}(x)^{-1}hxh^{-1}.  \label{eq:m(a,h)}
\end{equation}

\begin{theorem}
\label{thm:sigma_properties2}Let $G=BH$ be a transversal decomposition. For
all $x\in B$, $h\in H$, 
\begin{eqnarray}
m(x,1) &=&1  \label{eq:m(a,1)=1} \\
m(1,h) &=&1  \label{eq:m(1,h)=1} \\
\sigma _{m(x,h)} &=&\mu _{x}(\sigma _{h})  \label{eq:sigma_m(a,h)=mu}
\end{eqnarray}
\end{theorem}

\begin{proof}
Since $\sigma _{1}=I$, (\ref{eq:m(a,1)=1}) follows. By Theorem \ref
{thm:sigma_properties}(2), (\ref{eq:m(1,h)=1}) follows. Finally, for $x,y\in
B$, $h\in H$, we use (\ref{eq:sigma}), (\ref{eq:m(a,h)}), (\ref
{eq:B_operation}), and (\ref{eq:deviation}) to compute 
\begin{eqnarray*}
(\sigma _{m(x,h)}(y))H &=&m(x,h)yH \\
&=&\sigma _{h}(x)^{-1}hxh^{-1}yH \\
&=&\sigma _{h}(x)^{-1}hx\sigma _{h}^{-1}(y)H \\
&=&\sigma _{h}(x)^{-1}h(x\cdot \sigma _{h}^{-1}(y))H \\
&=&\sigma _{h}(x)^{-1}\sigma _{h}(x\cdot \sigma _{h}^{-1}(y))H \\
&=&(\sigma _{h}(x)^{-1}\cdot \sigma _{h}(x\cdot \sigma _{h}^{-1}(y)))H \\
&=&\left[ \left( L_{\sigma _{h}(x)}^{-1}\sigma _{h}L_{x}\sigma
_{h}^{-1}\right) (y)\right] H \\
&=&(\mu _{x}(\sigma _{h})(y))H
\end{eqnarray*}
This establishes (\ref{eq:sigma_m(a,h)=mu}).
\end{proof}

\begin{corollary}
\label{coro:sigma(H)_transassociant}$\sigma (H)$ is a transassociant of $%
(B,\cdot )$.
\end{corollary}

\begin{proof}
This follows from Theorem \ref{thm:sigma_properties}(2), (\ref
{eq:sigma_l(a,b)=L(a,b)}), and (\ref{eq:sigma_m(a,h)=mu}).
\end{proof}

We now consider the $B$- and $H$-components of a product of elements of $G$.

\begin{proposition}
\label{prop:factored_product2}For all $x,y\in B$, $h,k\in K$, 
\begin{equation}
xhyk=(x\cdot \sigma _{h}(y))l(x,\sigma _{h}(y))m(y,h)hk.  \label{eq:G_op}
\end{equation}
\end{proposition}

\begin{proof}
This is a direct computation.
\end{proof}

Comparison of Proposition \ref{prop:factored_product2} with Proposition \ref
{prop:factored_product} suggests the following.

\begin{definition}
\label{defn-internal-semi}Let $(B,\cdot )$ be the left loop induced by a
transversal decomposition $G=BH$. Then we say that $G$ is an \emph{internal
semidirect product} of $(B,\cdot )$ with $H$.
\end{definition}

\begin{remark}
\label{rem:special-cases}

\begin{enumerate}
\item  If $l:B\times B\rightarrow H$ is trivial, but $m:B\times H\rightarrow
H$ is nontrivial, then $(B,\cdot )$ is a subgroup of $G$, and $G=BH$ is the
internal version of Jajcay's ``rotary product'' \cite{jajcay} of subgroups.

\item  If $m:B\times H\rightarrow H$ is trivial, but $l:B\times B\rightarrow
H$ is nontrivial, then the product of $xh,yk\in G$ simplifies to 
\begin{equation*}
(xh)(yk)=(x\cdot \sigma _{h}(y))l(x,y)hk.
\end{equation*}
As will be shown below, in this case $(B,\cdot )$ is an $A_{l}$ left loop,
and $G=BH$ is the internal version of the semidirect product rediscovered
(for $A_{l}$, LIP left loops) by Kikkawa \cite{kikkawa} and Ungar \cite
{u-wags}.

\item  Both $l:B\times B\rightarrow H$ and $m:B\times H\rightarrow H$ are
trivial if and only if $B$ is a normal subgroup of $G$. In this case, $G=BH$
is the usual internal semidirect product of subgroups.

\item  Another case where $(B,\cdot )$ is a group is if $\sigma (H)=\{I\}$;
this follows from (\ref{eq:sigma_l(a,b)=L(a,b)}). However, if the
transversal mapping $l:B\times B\rightarrow H$ is nontrivial, then $B$ is
not a subgroup of $G$, and if $m:B\times H\rightarrow H$ is nontrivial, then 
$H$ does not normalize $B$.
\end{enumerate}
\end{remark}

We now consider the inheritance of internal semidirect product structure by
subgroups. Let $G=BH$ be a transversal decomposition giving an internal
semidirect product of $(B,\cdot )$ by $H$. Assume that $G_{1}$ is a subgroup
of $G$ and let $B_{1}=B\cap G_{1}$ and $H_{1}=H\cap G_{1}$. For $g\in G_{1}$%
, if $g=xh$ with $x\in B$, $h\in H$, then we clearly have $x\in B_{1}$ if
and only if $h\in H_{1}$. When either of these conditions hold, we say that $%
G_{1}$ \emph{respects} the transversal decomposition of $G$, or
equivalently, that $G_{1}$ respects the internal semidirect product
structure of $G$. If $G_{1}$ respects $G=BH$, then $G_{1}=B_{1}H_{1}$ is
itself a transversal decomposition, which means that $G_{1}$ is an internal
semidirect product of $B_{1}$ with $H_{1}$. In particular, the operation $%
\cdot $ on $B$ restricts to $B_{1}$, which shows that the left loop $%
(B_{1},\cdot )$ is a subloop of $(B,\cdot )$. Finally, if $G_{1}$ and $G_{2}$
are both subgroups respecting $G=BH$, then clearly the intersection $%
G_{1}\cap G_{2}$ satisfies this property as well.

Next we consider the inheritance of internal semidirect product structure by
factor groups. Let $G=BH$ be a transversal decomposition and let $%
K\vartriangleleft H$ be a normal subgroup of $G$. An arbitrary element $gK$
of $G/K$ factors as $gK=(xh)K=(xK)(hK)$ where $xK\in B_{K}=\{xK:x\in B\}$
and $hK\in H/K$. This factorization is clearly unique. Also, $B_{K}\cap
H/K=\{K\}$. Thus 
\begin{equation*}
G/K=B_{K}\ (H/K)
\end{equation*}
is a transversal decomposition of the factor group $G/K$. Denote the induced
binary operation (\ref{eq:B_operation}) by $\cdot _{K}:B_{K}\times
B_{K}\rightarrow B_{K}$ and the induced transversal mapping by $%
l_{K}:B_{K}\times B_{K}\rightarrow H/K$.

Since $B\cap K=\{1\}$, the set $B_{K}$ can be identified with $B$ itself.
Thus we compare two factorizations of products. For $x,y\in B$, we have 
\begin{equation*}
\left( xK\right) \left( yK\right) =\left( xK\cdot _{K}yK\right) l_{K}(xK,yK),
\end{equation*}
and also 
\begin{eqnarray*}
\left( xK\right) \left( yK\right) &=&(xy)K=\left( x\cdot y\right) l(x,y)K \\
&=&\left( \left( x\cdot y\right) K\right) \left( l(x,y)K\right) .
\end{eqnarray*}
By uniqueness, we have 
\begin{eqnarray*}
xK\cdot _{K}yK &=&\left( x\cdot y\right) K \\
l_{K}(xK,yK) &=&l(x,y)K
\end{eqnarray*}
for all $x,y\in B$. It follows that under the mapping $x\mapsto xK$, the
left loop $(B,\cdot )$ induced by the transversal decomposition $G=BH$ is
isomorphic to the left loop $(B_{K},\cdot _{K})$ induced by the transversal
decomposition $G/K=B_{K}\ (H/K)$. Making this identification, we may think
of 
\begin{equation}
G/K=B\ (H/K)  \label{eq:factor_decomp}
\end{equation}
as being a transversal decomposition of $G/K$.

As a specific example of factor group inheritance, let $G=BH$ be an internal
semidirect product of the left loop $(B,\cdot )$ with the subgroup $H$, and
let $K=\mathrm{core}_{G}(H)$. In addition to (\ref{eq:factor_decomp}), we
may make an additional observation. By Corollary \ref
{coro:sigma(H)_transassociant}, $\sigma (H)$ is a transassociant, and thus
we may form the standard semidirect product $B\rtimes \sigma (H)$. Thus the
exact sequence of groups 
\begin{equation}
1\rightarrow \mathrm{core}_{G}(H)\rightarrow H\rightarrow \sigma
(H)\rightarrow 1  \label{eq:H_seq}
\end{equation}
induces an exact sequence of semidirect product groups 
\begin{equation}
1\rightarrow \mathrm{core}_{G}(H)\rightarrow G\rightarrow B\rtimes \sigma
(H)\rightarrow 1.  \label{eq:internal_seq}
\end{equation}
The exactness of (\ref{eq:H_seq}) and (\ref{eq:internal_seq}) imply the
isomorphisms $H/\mathrm{core}_{G}(H)\cong \sigma (H)$ and $G/\mathrm{core}%
_{G}(H)\cong B\rtimes \sigma (H)$ Note also the obvious isomorphism of
groups $B\cdot H/\ker (\sigma )\rightarrow B\rtimes \sigma (H)$ given by $%
x(h\ker (\sigma ))\mapsto (x,\sigma _{h})$.

\begin{remark}
\label{rem:special-sequence}Consider again the special case of Remark \ref
{rem:special-cases}(4) where $\sigma (H)=\{I\}$, i.e., $H\vartriangleleft G$%
. Then making the usual identifications, (\ref{eq:internal_seq}) simplifies
to 
\begin{equation}
1\rightarrow H\rightarrow G\rightarrow B\rightarrow 1.
\label{eq:special_seq}
\end{equation}
As noted, if $l$ is nontrivial, then $B$ is a group, but not a subgroup.
Instead, we see from (\ref{eq:special_seq}) that $B$ is an isomorphic copy
of the factor group $G/H$.
\end{remark}

Let $G=BH$ be a transversal decomposition. Let $G_{0}=\langle B\rangle $,
and let $H_{0}=G_{0}\cap H$. It is straightforward to show that $%
H_{0}=\langle l(B,B)\rangle $ where $l(B,B)=\left\{ l(x,y):x,y\in B\right\} $%
. If $H$ is corefree, i.e., $\mathrm{core}_{G}(H)=\left\{ 1\right\} $, then
the restriction of $\sigma $ to $l(B,B)$ is a bijection onto the set $%
\left\{ L(x,y):x,y\in B\right\} $ of left inner mappings. Thus $\sigma
|_{H_{0}}$ is an isomorphism onto the left inner mapping group $\mathrm{LMlt}%
_{1}(B,\cdot )$. Putting these considerations together, we have the
following result \cite{phillips2}.

\begin{proposition}
\label{prop:G=LMlt}Let $G=BH$ be a transversal decomposition. Then $G\cong 
\mathrm{LMlt}(B,\cdot )$ if and only if $H$ is corefree and $G=\langle
B\rangle $.
\end{proposition}

There is a similar characterization of the multiplication group of a loop;
see Niemenmaa and Kepka \cite{nk}. In the more general case where the core
is not required to be trivial, Phillips has shown that a group $G=BH$ with $%
G=\langle B\rangle $ can be viewed as a left \emph{relative} multiplication
group of $(B,\cdot )$ as a subloop in some larger left loop \cite{phillips}.

We will consider some examples of internal semidirect products in the next
section, after using transversal decompositions to study certain varieties
of left loops.

\section{Varieties of Left Loops}

We begin by reviewing the definitions and properties of the varieties we
will consider. Let $(B,\cdot )$ be a left loop. $(B,\cdot )$ is said to
satisfy the \emph{left inverse property} (LIP) if 
\begin{equation}
L_{x}^{-1}=L_{x^{\rho }}  \label{eq:LIP}
\end{equation}
for all $x\in B$. This implies $x^{\rho }$ is a (unique) two-sided inverse
of $x\in B$. $(B,\cdot )$ is said to satisfy the \emph{left alternative
property} (LAP) if 
\begin{equation}
L_{x}L_{x}=L_{x\cdot x}  \label{eq:LAP}
\end{equation}
for all $x\in B$. $(B,\cdot )$ is said to be a (left) \emph{Bol loop} if 
\begin{equation}
L_{x}L_{y}L_{x}=L_{x\cdot \left( y\cdot x\right) }  \label{eq:Bol}
\end{equation}
for $x,y\in B$. A Bol loop satisfies LIP (take $x=y^{\rho }$ in (\ref{eq:Bol}%
)), LAP (take $y=1$ in (\ref{eq:Bol})), and is also a right loop \cite
{sharma}. If $B$ is a Bol loop, then for all $x,y\in B$, $L(x,y)$ is a
pseudo-automorphism with companion $(x\cdot y)\cdot (x^{\rho }\cdot y^{\rho
})$ \cite{good-robin}. Thus every Bol loop has the pseudo-$A_{l}$ property.
An $A_{l}$ Bol loop is sometimes called a ``gyrogroup'', as defined by Ungar
in \cite{u-axioms}. The equivalence of gyrogroups and $A_{l}$ Bol loops was
noted by R\`{o}zga \cite{rozga}.

There is an interesting intermediate variety of left loops which is defined
by the following identity: for all $x,y\in B$, 
\begin{equation}
L_{x}L_{y\cdot y}L_{x}=L_{x\cdot ((y\cdot y)\cdot x)}.  \label{eq:new_ident}
\end{equation}
Every such left loop has LAP (take $x=1$) and every Bol loop clearly
satisfies (\ref{eq:new_ident}). If $(B,\cdot )$ satisfies (\ref{eq:new_ident}%
) and the squaring mapping $x\longmapsto x^{2}$ is surjective, then $%
(B,\cdot )$ is a Bol loop.

A left loop $(B,\cdot )$ is said to satisfy the \emph{automorphic inverse
property} (AIP) if 
\begin{equation}
(x\cdot y)^{^{\rho }}=x^{\rho }\cdot y^{\rho }  \label{eq:AIP}
\end{equation}
for all $x,y\in B$, or simply $\rho L_{x}=L_{x^{\rho }}\rho $, where $\rho
:x\longmapsto x^{\rho }$. If $\rho $ is a permutation on $B$ with inverse $%
\lambda :x\longmapsto x^{\lambda }$, then AIP is characterized by $\rho
L_{x}\lambda \in L(B)$ for all $x\in B$. (Indeed, if $\rho L_{x}\lambda $ is
a translation, then applying it to $1$, we see that $\rho L_{x}\lambda
=L_{x^{\rho }}$.) An $A_{l}$, LIP, AIP left loop is called a \emph{Kikkawa
left loop} \cite{kiechle}.

An identity closely related to AIP is 
\begin{equation}
L_{x}L_{y}L_{y}L_{x}=L_{(x\cdot y)}L_{(x\cdot y)}  \label{eq:bruck1}
\end{equation}
for all $x,y\in B$. Applying both sides to $1$, we have 
\begin{equation}
x\cdot (y\cdot (y\cdot x)=(x\cdot y)\cdot (x\cdot y)  \label{eq:bruck2}
\end{equation}
for all $x,y\in B$. Taking $y=x^{\rho }$, we obtain $x^{\rho }=x^{\rho
}\cdot (x^{\rho }\cdot x)$. Cancelling, we have that every element of $B$
has a two-sided inverse.

The relationship between AIP, (\ref{eq:bruck1}) and our other identities is
summarized in the following result.

\begin{theorem}
\label{thm:AIP_equivs}Let $(B,\cdot )$ be a left loop.

\begin{enumerate}
\item  If $(B,\cdot )$ satisfies (\ref{eq:bruck1}), then LIP and AIP are
equivalent.

\item  If $(B,\cdot )$ is a Kikkawa left loop, then (\ref{eq:bruck1}) holds.

\item  If $(B,\cdot )$ satisfies LAP and (\ref{eq:bruck2}), then (\ref
{eq:bruck1}) and (\ref{eq:new_ident}) are equivalent.

\item  If $(B,\cdot )$ is a Bol loop, then AIP and (\ref{eq:bruck1}) are
equivalent.
\end{enumerate}
\end{theorem}

\begin{proof}
\noindent 1. Applying both sides of (\ref{eq:bruck1}) to $x\backslash
y^{\rho }$, we have $x\cdot y=(x\cdot y)\cdot ((x\cdot y)\cdot (x\backslash
y^{\rho }))$. Cancelling, we obtain $(x\cdot y)^{\rho }=x\backslash y^{\rho
} $ for all $x,y\in B$, or equivalently, $\rho L_{x}=L_{x}^{-1}\rho $ for
all $x\in B$. If LIP\ holds, then $\rho L_{x}=L_{x^{\rho }}\rho $ for all $%
x\in B$, which is AIP. Conversely, if AIP\ holds, then $L_{x^{\rho }}\rho
=\rho L_{x}=L_{x}^{-1}\rho $. Cancelling $\rho $, we obtain LIP.

\noindent 2. See \cite{kikkawa}, Prop. 1.13.

\noindent 3. Using LAP, (\ref{eq:bruck2}), and LAP again, we have 
\begin{equation*}
L_{x\cdot y}L_{x\cdot y}=L_{(x\cdot y)\cdot (x\cdot y)}=L_{x\cdot (y\cdot
(y\cdot x)}=L_{x\cdot ((y\cdot y)\cdot x)}
\end{equation*}
for all $x,y\in B$. From this the equivalence is clear.

\noindent 4. From the remarks following (\ref{eq:Bol}), we have that in a
Bol loop satisfying AIP, each left inner mapping $L(x,y)$ is an
automorphism, i.e., the $A_{l}$ property holds. By (2), (\ref{eq:bruck1})
holds. The converse follows from (1), since every Bol loop has LIP.
\end{proof}

For a related discussion, see also \cite{kiechle}, especially pp. 37-38 and
58-59. The implication ``AIP implies LIP'' in Theorem \ref{thm:AIP_equivs}%
(1) seems to be new, as is the connection with the variety of left loops
satisfying (\ref{eq:new_ident}).

A Bol loop satisfying AIP or (\ref{eq:bruck1})\ is called a \emph{Bruck loop}%
. This was a term coined by Robinson in his dissertation \cite{rob-diss}
under the additional assumption that the mapping $x\longmapsto x\cdot x$ is
a permutation. The term acquired its generally accepted present meaning due
to a remark of Glauberman (\cite{glauberman}, p.376), who also coined the
term \emph{B-loop} to describe finite, odd order, Bruck loops. The term
``B-loop'' is now used for the general (not necessarily finite) case in
which squaring is a permutation. Bruck loops are also known as ``K-loops'' 
\cite{kar-wef} \cite{kiechle} \cite{kk} \cite{kreuzer} \cite{u-wags} and as
``(gyrocommutative) gyrogroups'' \cite{u-grouplike} \cite{u-axioms}. The
direct equivalence between Bruck loops and K-loops was shown by Kreuzer \cite
{kreuzer}. The direct equivalence between Bruck loops and gyrocommutative
gyrogroups was shown by Sabinin \textit{et al} \cite{sss}, and was also
noted by R\`{o}zga \cite{rozga}. (The direct equivalence of K-loops and
gyrocommutative gyrogroups was a well-known folk result.)

Let $G=BH$ be a transversal decomposition. In the following discussion we
will abbreviate the core of $H$ in $G$ by $N=\mathrm{core}_{G}(H)$. For $%
B\subseteq G$, let $B^{2}=\left\{ x^{2}:x\in B\right\} $ and let $%
B^{-1}=\left\{ x^{-1}:x\in B\right\} $. We introduce the following
conditions on the transversal $B$.

\begin{center}
\begin{tabular}{ll}
(G-LIP) & $B^{-1}\subseteq BN.$ \\ 
(G-LAP) & $B^{2}\subseteq BN.$ \\ 
(G-Bol) & For all $x\in B$, $xBx\subseteq BN.$ \\ 
(G-W) & For all $x\in B$, $xB^{2}x\subseteq BN.$ \\ 
(G-Br) & For all $x,y\in B$, $xy^{2}x\in (x\cdot y)^{2}N.$%
\end{tabular}
\end{center}

In the following result, we show that these transversal conditions
characterize the aforementioned left loop varieties. In this respect, our
approach is sharper than that of other authors \cite{u-wags} \cite
{kreuz-wefel}. For instance, it is well-known that the property $%
B^{-1}\subseteq B$ implies that the left loop $(B,\cdot )$ has LIP. However,
this characterizes LIP only in decompositions $G=BH$ in which $H$ is
corefree.

Recall the notation $B_{N}=\left\{ xN:x\in B\right\} $ and that the left
loop $(B_{N},\cdot _{N})$ is an isomorphic copy of $(B,\cdot )$.

\begin{theorem}
\label{thm:internal_idents}Let $G=BH$ be a transversal decomposition. In
each of the following, assertions (i), (ii), and (iii) are equivalent.

\begin{enumerate}
\item  (i) (G-LIP); (ii) (G/N-LIP); (iii) $(B,\cdot )$ satisfies LIP.

\item  (i) (G-LAP); (ii) (G/N-LAP); (iii) $(B,\cdot )$ satisfies LAP.

\item  (i) (G-Bol); (ii) (G/N-Bol); (iii) $(B,\cdot )$ is a Bol loop.

\item  (i) (G-W); (ii) (G/N-W); (iii) $(B,\cdot )$ satisfies (\ref
{eq:new_ident}).

\item  (i) (G-Br); (ii) (G/N-Br); (iii) $(B,\cdot )$ satisfies (\ref
{eq:bruck1}).
\end{enumerate}
\end{theorem}

\begin{proof}
1. For all $x\in B$, we have $L_{x}L_{x^{\rho }}=\sigma _{x}\sigma _{x^{\rho
}}=\sigma _{xx^{\rho }}$. By (\ref{eq:LIP}), $(B,\cdot )$ has LIP if and
only if, for each $x\in B$, $\sigma _{xx^{\rho }}=I$, or equivalently, $%
xx^{\rho }\in N$. This is equivalent to $x^{-1}\in x^{\rho }N$. Since $%
x^{-1}H=x^{\rho }H$, LIP holds if and only if, for each $x\in B$, $x^{-1}\in
BN$, which is (G-LIP), or equivalently, $x^{-1}N\in B_{N}$, which is
(G/N-LIP).

\noindent 2. For all $x\in B$, we have $L_{(x\cdot x)}^{-1}L_{x}L_{x}=\sigma
_{(x\cdot x)^{-1}x^{2}}$. By (\ref{eq:LAP}), $(B,\cdot )$ has LAP if and
only if, for each $x\in B$, $\sigma _{(x\cdot x)^{-1}x^{2}}=I$, or
equivalently, $(x\cdot x)^{-1}x^{2}\in N$, or $x^{2}\in (x\cdot x)N$. Since $%
x^{2}H=(x\cdot x)H$, LAP holds if and only if, for each $x\in B$, $x^{2}\in
BN$, which is (G-LAP), or equivalently, $x^{2}N\in B_{N}$, which is
(G/N-LAP).

\noindent 3. For all $x,y\in B$, we have $L_{x\cdot (y\cdot
x)}^{-1}L_{x}L_{y}L_{x}=\sigma _{(x\cdot (y\cdot x))^{-1}xyx}$. By (\ref
{eq:Bol}), $(B,\cdot )$ is a Bol loop if and only if, for every $x,y\in B$, $%
\sigma _{(x\cdot (y\cdot x))^{-1}xyx}=I$, or equivalently, $(x\cdot (y\cdot
x))^{-1}xyx\in N$, or $xyx\in (x\cdot (y\cdot x))N$. Since $xyxH=(x\cdot
(y\cdot x))H$, $(B,\cdot )$ is a Bol loop if and only if, for every $x,y\in
B $, $xyx\in BN$, which is (G-Bol), or equivalently, for every $x,y\in B$, $%
(xN)(yN)(xN)=x(yN)x\in B_{N}$, which is (G/N-Bol).

\noindent 4. The proof is similar to that of (3), \textit{mutatis mutandis}.

\noindent 5. For all $x,y\in B$, we have $L_{(x\cdot y)}^{-1}L_{(x\cdot
y)}^{-1}L_{x}L_{y}L_{y}L_{x}=\sigma _{(x\cdot y)^{-2}xy^{2}x}$. Thus $%
(B,\cdot )$ satisfies (\ref{eq:bruck1}) if and only if, for every $x,y\in B$%
, $(x\cdot y)^{-2}xy^{2}x\in N$, or $xy^{2}x\in (x\cdot y)^{2}N$, which is
(G-Br), or equivalently, $(xN)(yN)^{2}(xN)=(xN\cdot _{N}yN)^{2}$, which is
(G/N-Br).
\end{proof}

\begin{corollary}
\label{coro-loop-identities}Let $(B,\cdot )$ be a left loop.

\begin{enumerate}
\item  $(B,\cdot )$ satisfies LIP if and only if, for all $x\in B$, $%
L_{x}^{-1}\in L(B)$.

\item  $(B,\cdot )$ satisfies LAP if and only if, for all $x\in B$, $%
L_{x}^{2}\in L(B)$.

\item  $(B,\cdot )$ is a Bol loop if and only if, for all $x,b\in B$, $%
L_{x}L_{y}L_{x}\in L(B)$.

\item  $(B,\cdot )$ is satisfies (\ref{eq:new_ident}) if and only if, for
all $x,y\in B$, $L_{x}L_{y\cdot y}L_{x}\in L(B)$.
\end{enumerate}
\end{corollary}

\begin{remark}
\label{rem-extensions}

\begin{enumerate}
\item  Other varieties of left loops may be characterized in terms of
internal semidirect product structure. For instance, a left loop $(B,\cdot )$
said to have the \emph{LC property} if $L_{x}L_{x}L_{y}=L_{x\cdot (x\cdot
y)} $ for all $x,y\in B$ \cite{fenyves}. If $G=BH$ is a transversal
decomposition and $N=\mathrm{core}_{G}(H)$, then $(B,\cdot )$ is an LC left
loop if and only if $x^{2}B\subseteq BN$ for all $x\in B$. Similarly, a left
loop $(B,\cdot )$ is said to be \emph{left conjugacy closed} if $%
L_{x}L_{y}L_{x}^{-1}\in L(B)$ for all $x,y\in B$ \cite{nag-stram}. If $G=BH$
is a transversal decomposition, then $(B,\cdot )$ is left conjugacy closed
if and only if $xBx^{-1}\subseteq BN$ for all $x\in B$. The proofs of both
of these assertions are similar to that of Theorem \ref{thm:internal_idents}%
(3), \textit{mutatis mutandis}.

\item  Let $G$ be a group. A subset $B\subseteq G$ is said to be a \emph{%
twisted subgroup} of $G$ if $1\in B$ and if $xBx\subseteq B$ for all $x\in B$
\cite{aschbacher} \cite{foguel-ungar}. In this jargon, we can restate
Theorem \ref{thm:internal_idents}(3) as follows: \emph{If} $G=BH$ \emph{is a
transversal decomposition with }$H$ \emph{corefree, then} $B$ \emph{is a
twisted subgroup if and only if }$(B,\cdot )$ \emph{is a Bol loop}. This
generalizes \cite{foguel-ungar}, Thm. 3.8.

\item  Let $G=BH$ with $H$ corefree. Assume $G$ satisfies (G-Bol) and
(G-Br). By Theorem \ref{thm:internal_idents}(3),(5), $(B,\cdot )$ is a Bruck
loop. If the squaring mapping $B\rightarrow B:x\mapsto x^{2}$ is a
permutation, then $(B,\cdot )$ is a B-loop, and the loop operation in $B$ is
given by 
\begin{equation}
x\cdot y=(xy^{2}x)^{1/2},  \label{eq:a.b=sqrt(abba)}
\end{equation}
$x,y\in B$. It is interesting to compare the operation (\ref
{eq:a.b=sqrt(abba)}) in $B$ with the operation discussed by Glauberman \cite
{glauberman} and Kikkawa \cite{kikkawa} (Ex. 1.5), namely, $x\ast
y=x^{1/2}yx^{1/2}$. Observe that $x\cdot y=(x^{2}\ast y^{2})^{1/2}$. The
operation $\cdot $ of the present paper is more natural with respect to the
group structure $G=BH$. See also \cite{kiechle}, pp. 60-62.
\end{enumerate}
\end{remark}

The significance of Theorem \ref{thm:internal_idents} is that in the study
of these particular varieties of left loops via internal semidirect
products, there is no loss in assuming that the subgroup $H$ is corefree. In
particular, for simply proving facts about abstract left loops, one can
always work with the transversal $L(B)$ in the left multiplication group $%
\mathrm{LMlt}(B,\cdot )$. On the other hand, for examples it is frequently
the case that the group $G$ and the transversal decomposition $G=BH$ are
given. In this instance it is preferable to work directly with the group $G$.

Let $G=BH$ be a transversal decomposition, and again let $N=\mathrm{core}%
_{G}(H)$. We now consider conditions on the interaction between $B$ and $H$.

\begin{center}
\begin{tabular}{ll}
(G-PsA$_{l}$) & For each $h\in H$, $chBh^{-1}\subseteq BN$ for some $c\in B$.
\\ 
(G-A$_{l}$) & For all $h\in H$, $hBh^{-1}\subseteq BN$.
\end{tabular}
\end{center}

\begin{lemma}
\label{lem:pseudo_lemma}Let $G=BH$ be a transversal decomposition, and let $%
h\in H$ be given. In each of the following, (i), (ii), and (iii) are
equivalent.

\begin{enumerate}
\item  (i) $\sigma _{h}\in \mathrm{PsAut}(B,\cdot )$; (ii) $%
chBh^{-1}\subseteq BN$ for some $c\in B$; (iii) for some $c\in B$, $%
l(c,\sigma _{h}(x))m(x,h)\in N$ for all $x\in B$.

\item  (i) $\sigma _{h}\in \mathrm{Aut}(B,\cdot )$; (ii) $hBh^{-1}\subseteq
BN$; (iii) $m(x,h)\in N$ for all $x\in B$.
\end{enumerate}
\end{lemma}

\begin{proof}
\noindent 1. For $x,c\in B$, we have 
\begin{equation*}
chxh^{-1}=c\sigma _{h}(x)m(x,h)=(c\cdot \sigma _{h}(x))l(c,\sigma
_{h}(x))m(x,h).
\end{equation*}
Thus the equivalence of (ii) and (iii) is clear. For $x,y,c\in B$, we have 
\begin{equation*}
(c\cdot \sigma _{h}(x\cdot y))H=chxyH=chxh^{-1}hyH=chxh^{-1}\sigma _{h}(y)H,
\end{equation*}
and 
\begin{equation*}
((c\cdot \sigma _{h}(x))\cdot \sigma _{h}(y))H=(c\cdot \sigma _{h}(x))\sigma
_{h}(y)H.
\end{equation*}
Thus $\sigma _{h}$ is a pseudo-automorphism with companion $c$ if and only
if, for all $x\in B$, $(c\cdot \sigma _{h}(x))^{-1}chxh^{-1}\in N$, or
equivalently, $chxh^{-1}\in (c\cdot \sigma _{h}(x))N$. This establishes the
equivalence of (i) and (ii).

\noindent 2. This follows from (1) with $c=1$ (using (\ref
{eq:l(1,a)=l(a,1)=1}) for (iii)).
\end{proof}

Recall the notation $G_{0}=\langle B\rangle $.

\begin{theorem}
\label{thm:aut_char}Let $G=BH$ be a transversal decomposition.

\begin{enumerate}
\item  (G-PsA$_{l}$) holds if and only if (G/N-PsA$_{l}$) holds.

\item  (G$_{0}$-PsA$_{l}$) holds if and only if $(B,\cdot )$ is a pseudo-A$%
_{l}$ left loop.

\item  (G-A$_{l}$) holds if and only if (G/N-A$_{l}$) holds.

\item  (G$_{0}$-A$_{l}$) holds if and only if $(B,\cdot )$ is an A$_{l}$
left loop.
\end{enumerate}
\end{theorem}

\begin{proof}
\noindent 1. (G-PsA$_{l}$) holds if and only if, for every $x\in B$, $h\in H$%
, there exists $c\in B$ such that $chxh^{-1}\in (c\cdot \sigma _{h}(x))N$,
or 
\begin{equation*}
(cN)(hN)(xN)(h^{-1}N)=chxh^{-1}N=(c\cdot \sigma _{h}(x))N\in B_{N}.
\end{equation*}
This is (G/N-PsA$_{l}$).

\noindent 2. This follows from applying Lemma \ref{lem:pseudo_lemma}(1)\ to $%
G_{0}$.

\noindent 3. The proof is similar to that of (1), \textit{mutatis mutandis}.

\noindent 4. This follows from applying Lemma \ref{lem:pseudo_lemma}(2)\ to $%
G_{0}$.
\end{proof}

Let $G=BH$ be a transversal decomposition. We define a mapping on $G$ as
follows: 
\begin{equation}
\tau :G\rightarrow G:ah\longmapsto a^{-1}h  \label{eq:tau}
\end{equation}
for $a\in B$, $h\in H$. Recall that in $(B,\cdot )$, we also have the right
inversion mapping $\rho :B\rightarrow B:x\longmapsto x^{\rho }$. The
following result shows that $\tau $ encodes information about $\rho $ at the
level of the entire group $G$.

\begin{theorem}
\label{thm:tau_props}Let $G=BH$ be a transversal decomposition.

\begin{enumerate}
\item  $\tau $ is injective if and only if $\rho $ is injective.

\item  $\tau $ is surjective if and only if $\rho $ is surjective, i.e., if
and only if every element of $(B,\cdot )$ has a left inverse.

\item  $\tau $ is bijective if and only if $\rho $ is bijective, i.e., if
and only if every element of $(B,\cdot )$ has a unique left inverse.

\item  $\tau ^{2}=I$ if and only if, for all $x\in B$, $xx^{\rho }x^{\rho
}x=e$. If these conditions hold, then every element of $(B,\cdot )$ has a
two-sided inverse.
\end{enumerate}
\end{theorem}

\begin{proof}
\noindent 1. Assume $\rho $ is injective and $\tau (xh)=\tau (yk)$ for some $%
x,y\in B$, $h,k\in H$. Then $x^{\rho }l(x,x^{\rho })^{-1}h=y^{\rho
}l(y,y^{\rho })^{-1}k$, using (\ref{eq:l(a,ar)=aar}). Matching $B$ and $H$
components, this is equivalent to $x^{\rho }=y^{\rho }$ and\ $l(x,x^{\rho
})^{-1}h=l(y,y^{\rho })^{-1}k$. Now $x^{\rho }=y^{\rho }$ implies $x=y$, and
cancelling $l(x,x^{\rho })^{-1}=l(y,y^{\rho })^{-1}$, we have $h=k$. Thus $%
\tau $ is injective. Conversely, assume $\tau $ is injective and $x^{\rho
}=y^{\rho }$. Then using (\ref{eq:l(a,ar)=aar}) again, 
\begin{equation*}
\tau (xl(x,x^{\rho }))=x^{-1}l(x,x^{\rho })=x^{\rho }=y^{\rho
}=y^{-1}l(y,y^{\rho })=\tau (yl(y,y^{\rho })).
\end{equation*}
Thus $xl(x,x^{\rho })=yl(y,y^{\rho })$. Matching components, $x=y$.

\noindent 2. $\tau $ is surjective if and only if, for each $a\in B$, there
exists $x\in B$ and $h\in H$ such that $a=\tau (xh)=x^{-1}h$, i.e., $xa=h\in
H$. This implies $x\cdot a=e$, and conversely, if $x\cdot a=e$, then $xa=h$
for some $h\in H$. This establishes the desired equivalence.

\noindent 3. This follows from (1) and (2).

\noindent 4. Using (\ref{eq:l(a,ar)=aar}), we have for all $x\in B$, 
\begin{equation*}
\tau ^{2}(x)=\tau (x^{-1})=\tau (x^{\rho }l(x,x^{\rho })^{-1})=x^{\rho \rho
}l(x^{\rho },x^{\rho \rho })^{-1}l(x,x^{\rho })^{-1}.
\end{equation*}
Now 
\begin{equation*}
l(x^{\rho },x^{\rho \rho })^{-1}l(x,x^{\rho })^{-1}=(l(x,x^{\rho })l(x^{\rho
},x^{\rho \rho }))^{-1}=(xx^{\rho }x^{\rho }x^{\rho \rho })^{-1}.
\end{equation*}
If $\tau ^{2}(x)=x$, then matching components gives $x=x^{\rho \rho }$
(i.e., $x^{\rho }$ is a two-sided inverse of $x$), and $xx^{\rho }x^{\rho
}x^{\rho \rho }=xx^{\rho }x^{\rho }x=e$. Conversely, if $xx^{\rho }x^{\rho
}x=e$ for all $x\in B$, then $\tau ^{2}(x)=x^{\rho \rho }([xx^{\rho }x^{\rho
}]x^{\rho \rho })^{-1}=x^{\rho \rho }(x^{-1}x^{\rho \rho })^{-1}=x$.
\end{proof}

\begin{remark}
\label{rem:tau_lip}If $G=BH$ is a transversal decomposition, then obviously
(G-LIP) holds if and only if $\tau (B)\subseteq BN$. In previous studies of
the mapping $\tau $ (e.g., \cite{kar-wef}, \cite{kikkawa}), (G-LIP) was
assumed at the outset. Theorem \ref{thm:tau_props} shows that $\tau $
contains information about $(B,\cdot )$ even under much weaker conditions.
\end{remark}

Now assume that $\tau :G\rightarrow G$ is a semi-automorphism, i.e., 
\begin{equation}
\tau (g_{1}g_{2}g_{1})=\tau (g_{1})\tau (g_{2})\tau (g_{1})
\label{eq:tau_semi}
\end{equation}
for all $g_{1},g_{2}\in G$. For $x\in B$, we have 
\begin{equation*}
x^{-1}=\tau (x)=\tau (xx^{-1}x)=\tau (x)\tau (x^{-1})\tau (x)=x^{-1}\tau
(x^{-1})x^{-1}.
\end{equation*}
Cancelling, we obtain $\tau ^{2}(x)=\tau (x^{-1})=x$. By Theorem \ref
{thm:tau_props}(4), each element of $(B,\cdot )$ has a two-sided inverse.
Now for $x,y\in B$, we compute 
\begin{eqnarray*}
x^{-1}y^{-1}x^{-1} &=&\tau (x)\tau (y)\tau (x)=\tau (xyx)=\tau ((x\cdot
(y\cdot x))l(x,y\cdot x)l(y,x)) \\
&=&(x\cdot (y\cdot x))^{-1}l(x,y\cdot x)l(y,x)=(x\cdot (y\cdot x))^{-2}xyx.
\end{eqnarray*}
Thus 
\begin{equation}
(xyx)^{2}=(x\cdot (y\cdot x))^{2}.  \label{eq:squares}
\end{equation}

\begin{proposition}
\label{prop:tau_semi_bol}Let $G$ be such that the squaring mapping $%
g\longmapsto g^{2}$ is injective. If $\tau :G\rightarrow G$ is a
semi-automorphism, then $(B,\cdot )$ is a Bol loop.
\end{proposition}

\begin{proof}
From (\ref{eq:squares}), we have that $xyx=x\cdot (y\cdot x)\in B$ for all $%
x,y\in B$. This implies (G-Bol), and the result follows from Theorem \ref
{thm:internal_idents}(3).
\end{proof}

\begin{remark}
\label{rem:tau_trivial}Let $G=BH$ be a transversal decomposition with $N=%
\mathrm{core}_{G}(H)$. Clearly (G-Bol) and (G-A$_{l}$) hold if and only if $%
gB\tau (g)^{-1}\subseteq BN$ for all $g\in G$.
\end{remark}

Finally, we characterize those groups for which $\tau $ is an automorphism.

\begin{theorem}
\label{thm:tau_aut}Let $G=BH$ be a transversal decomposition.

\begin{enumerate}
\item  $\tau :G\rightarrow G$ is an automorphism, then for all $x,y\in B$, 
\begin{equation}
xy^{2}x=(x\cdot y)^{2}  \label{eq:strong_br}
\end{equation}
and for all $y\in B$, $h\in H$, 
\begin{equation}
\sigma _{h}(y)^{2}=hy^{2}h^{-1}.  \label{eq:sig_square}
\end{equation}

\item  If $H$ is corefree and (G-A$_{l}$) holds, then $\tau $ is an
automorphism if and only if (\ref{eq:strong_br}) holds.

\item  If $G=\langle B\rangle $, then $\tau $ is an automorphism if and only
if (\ref{eq:strong_br}) holds.
\end{enumerate}
\end{theorem}

\begin{proof}
\noindent 1. For $x,y\in B$, $h\in H$, we compute 
\begin{equation*}
\tau (xhyk)=\tau ((x\cdot \sigma _{h}(y))[(x\cdot \sigma
_{h}(y))^{-1}xhyk])=(x\cdot \sigma _{h}(y))^{-2}xhyk
\end{equation*}
and 
\begin{equation*}
\tau (xh)\tau (yk)=x^{-1}hy^{-1}k.
\end{equation*}
Thus $\tau $ is an automorphism if and only if 
\begin{equation}
(x\cdot \sigma _{h}(y))^{2}=xhy^{2}h^{-1}x  \label{eq:stronga}
\end{equation}
for all $x,y\in B$, $h\in H$. Taking $h=1$ gives (\ref{eq:strong_br}).
Taking $x=1$ gives (\ref{eq:sig_square}). Conversely, (\ref{eq:strong_br})
and (\ref{eq:sig_square}) clearly imply (\ref{eq:stronga}).

\noindent 2. If $H$ is corefree and (G-A$_{l}$) holds, then (\ref
{eq:sig_square}) holds, so the result follows from (1).

\noindent 3. Assume $G=\langle B\rangle $. We will show that (\ref
{eq:strong_br}) implies (\ref{eq:sig_square}). First note that (\ref
{eq:strong_br}) implies $y^{2}=(x\cdot (x\backslash y))^{2}=x(x\backslash
y)^{2}x$, or 
\begin{equation}
x^{-1}y^{2}x^{-1}=(x\backslash y)^{2}  \label{eq:strong_br2}
\end{equation}
for all $x,y\in B$. Fix $h\in H$. Then $h=a_{1}^{\varepsilon
_{1}}a_{2}^{\varepsilon _{2}}\cdots a_{n}^{\varepsilon _{n}}$ for some $%
a_{i}\in B$, where $\varepsilon _{i}=\pm 1$, $i=1,\ldots ,n$. Using (\ref
{eq:strong_br}) and (\ref{eq:strong_br2}), we have for $y\in B$, 
\begin{equation*}
hy^{2}h^{-1}=a_{1}^{\varepsilon _{1}}\cdots a_{n}^{\varepsilon
_{n}}y^{2}a_{n}^{-\varepsilon _{n}}\cdots a_{1}^{-\varepsilon
_{1}}=(a_{1}\cdot _{\varepsilon _{1}}\cdots (a_{n}\cdot _{\varepsilon
_{n}}y)\cdots )
\end{equation*}
where $\cdot _{\varepsilon _{i}}=\cdot $ if $\varepsilon _{i}=1$ and $\cdot
_{\varepsilon _{i}}=\backslash $ if $\varepsilon _{i}=-1$. Thus 
\begin{equation*}
hy^{2}h^{-1}=\sigma _{a_{1}}^{\varepsilon _{1}}\cdots \sigma
_{a_{n}}^{\varepsilon _{n}}(y^{2})=\sigma _{a_{1}^{\varepsilon _{1}}}\cdots
\sigma _{a_{n}^{\varepsilon _{n}}}(y^{2})=\sigma _{h}(y^{2}).
\end{equation*}
This completes the proof.
\end{proof}

Note that (\ref{eq:strong_br}) implies (G-Br), and reduces to (G-Br) if $H$
is corefree.

\begin{corollary}
\label{coro:tau_aut}Let $(B,\cdot )$ be a left loop, let $G=\mathrm{LMlt}%
(B,\cdot )$, and define $\tau :G\rightarrow G$ by $\tau (L_{x}\phi
)=L_{x}^{-1}\phi $ ($x\in B$, $\phi \in \mathrm{LMlt}_{1}(B,\cdot )$).

\begin{enumerate}
\item  $(B,\cdot )$ satisfies (\ref{eq:bruck1}) if and only $\tau $ is an
automorphism.

\item  If the squaring mapping $G\rightarrow G:g\longmapsto g^{2}$ is
injective, then $(B,\cdot )$ is a Bruck loop if and only if $\tau $ is an
automorphism.
\end{enumerate}
\end{corollary}

\begin{proof}
\noindent 1. This follows from Proposition \ref{prop:G=LMlt}, Theorem \ref
{thm:tau_aut}(2), and Theorem \ref{thm:internal_idents}(5).

\noindent 2. This follows from (1) and Proposition \ref{prop:tau_semi_bol}.
\end{proof}

\begin{remark}
Automorphisms (or, equivalently, anti-automorphisms of order $2$) have been
used by various authors to study Bruck loops. See, for instance, \cite
{foguel-ungar}, \cite{glauberman}, \cite{im}, \cite{kar-wef}, \cite{kikkawa}%
. Theorem \ref{thm:tau_aut} and its corollary clarify the exact relationship
between $\tau $ and (\ref{eq:bruck1}).
\end{remark}

We now consider some examples of transversal decompositions giving internal
semidirect products.

\begin{example}
\label{ex-polar}(\emph{Polar Decomposition}) Let $GL(n,\mathbb{C})$ denote
the general linear group of $n\times n$ complex invertible matrices, let $%
\mathcal{P}(n)$ denote the subset of all $n\times n$ positive definite
Hermitian matrices, and let $U(n)$ denote the subgroup of unitary $n\times n$
complex matrices. The \emph{polar decomposition} asserts that every $M\in
GL(n,\mathbb{C})$ can be uniquely factored as $M=AU$ for a unique $%
A=(MM^{\ast })^{1/2}\in \mathcal{P}(n)$ and $U=(MM^{\ast })^{-1/2}A\in U(n)$%
, where $M^{\ast }$ is the conjugate transpose of $M$ and where the unique
positive definite square root of $MM^{\ast }$ is intended. Thus the polar
decomposition is a transversal decomposition 
\begin{equation}
GL(n,\mathbb{C})=\mathcal{P}(n)\cdot U(n).  \label{eq:polar}
\end{equation}
The induced binary operation (\ref{eq:B_operation}), denoted here by $\odot $%
, is given by 
\begin{equation}
A\odot B=(AB(AB)^{\ast })^{1/2}=(AB^{2}A)^{1/2}  \label{eq:polar_op}
\end{equation}
for $A,B\in \mathcal{P}(n)$; compare with \ref{eq:a.b=sqrt(abba)}. The
transversal mapping $l:\mathcal{P}(n)\times \mathcal{P}(n)\rightarrow U(n)$
is given by 
\begin{equation}
l(A,B)=(AB^{2}A)^{-1/2}AB  \label{eq:polar_l}
\end{equation}
for $A,B\in \mathcal{P}(n)$. We have $\mathcal{P}(n)^{-1}\subseteq \mathcal{P%
}(n)$, i.e., (G-LIP) holds. The involution $\tau :GL(n,\mathbb{C}%
)\rightarrow GL(n,\mathbb{C}):AU\longmapsto A^{-1}U$ is given explicitly by 
\begin{eqnarray*}
\tau (M) &=&\tau ((MM^{\ast })^{1/2}(MM^{\ast })^{-1/2}M) \\
&=&(MM^{\ast })^{-1/2}(MM^{\ast })^{-1/2}M \\
&=&(MM^{\ast })^{-1}M \\
&=&(M^{\ast })^{-1}
\end{eqnarray*}
for $M\in GL(n,\mathbb{C})$. For $A\in \mathcal{P}(n)$, $M\in GL(n,\mathbb{C}%
)$, we have $MA\tau (M)^{-1}=MAM^{\ast }\in \mathcal{P}(n)$. By Remark \ref
{rem:tau_trivial}, (G-Bol) and (G-A$_{l}$) hold. Since $\tau $ is the
composition of two involutory anti-automorphisms (conjugate transposition
and inversion), $\tau $ is an automorphism. By Theorem \ref{thm:tau_aut}(3),
(G-Br) holds. In addition, the squaring mapping is a permutation of $%
\mathcal{P}(n)$. Putting all these facts together, we have that $(\mathcal{P}%
(n),\oplus )$ is a $B$-loop.
\end{example}

\begin{example}
\label{ex-subgroups}(\emph{Subgroups of }$GL(n,\mathbb{C})$) Subgroups of $%
GL(n,\mathbb{C})$ respecting the polar decomposition include the special
linear group $SL(n,\mathbb{C})$, the real general linear group $GL(n,\mathbb{%
R})$, the group $U(m,n)$, and the complex symplectic group $Sp(n,\mathbb{C})$%
. Taking intersections of these yields more such groups. Here we limit
ourselves to one specific example: The group $SU(1,1)$ consists of those $%
2\times 2$ complex matrices preserving the form $|z_{1}|^{2}-|z_{2}|^{2}$ on 
$\mathbb{C}^{2}$ which also have determinant $1$. The polar decomposition of
this group is 
\begin{equation*}
SU(1,1)=\mathcal{P}U(1,1)\cdot S(U(1)\times U(1)).
\end{equation*}
where $S(U(1)\times U(1))$ is the subgroup of matrices of the form $\left( 
\begin{array}{cc}
\bar{a} & 0 \\ 
0 & a
\end{array}
\right) $, $a\in S^{1}$, and $\mathcal{P}U(1,1)=\mathcal{P}(n)\cap U(1,1)$
is the set of positive definite Hermitian matrices in $U(1,1)$ (such
matrices necessarily have determinant $1$). Thus $(\mathcal{P}U(1,1),\oplus
) $ is a subloop of the $B$-loop $(\mathcal{P}(2),\oplus )$ of Example \ref
{ex-polar}. Note that the mapping $Q:S^{1}\rightarrow S(U(1)\times U(1))$
defined by $Q(a)=\left( 
\begin{array}{cc}
\bar{a} & 0 \\ 
0 & a
\end{array}
\right) \ $is an isomorphism of groups. A matrix $L$ in $\mathcal{P}U(1,1)$
can be parametrized by a number $z\in \mathbb{D}$ as follows: 
\begin{equation}
L=L(z)=\gamma _{z}\left( 
\begin{array}{cc}
1 & \bar{z} \\ 
z & 1
\end{array}
\right)  \label{eq:L(z)}
\end{equation}
where $\gamma _{z}=(1-|z|^{2})^{-1/2}$. The mapping $L:\mathbb{D}\rightarrow 
\mathcal{P}U(1,1)$ turns out to be an isomorphism from the $B$-loop $(%
\mathbb{D},\oplus )$ of Example \ref{ex-disk} to $(\mathcal{P}U(1,1),\oplus
) $.
\end{example}

\begin{remark}
\label{rem-cartan}The polar decompositions of Examples \ref{ex-polar} and 
\ref{ex-subgroups} are special cases of the global Cartan decomposition of a
Lie group associated with a Riemannian symmetric space of noncompact type 
\cite{helg}. Any such space, realized as a subset of the Lie group, can be
given the structure of a $B$-loop; this actually follows quite easily from
the results herein. For the Hermitian case (bounded symmetric domains), the
result was shown in \cite{fried-ungar}, while the general Riemannian case
was worked out in \cite{kram-urb} and \cite{rozga}. For related work over
Pythagorean fields, see \cite{kk}.
\end{remark}

\begin{example}
\label{ex-nongyrocomm}We will use the polar decomposition of the group $%
U(1,1)$ to construct a different transversal decomposition, leading to a
different internal semidirect product structure. The polar decomposition of
a given $A\in U(1,1)$ is $A=L(z)V(a,b)$, with $L(z)\in \mathcal{P}U(1,1)$
given by (\ref{eq:L(z)}) for some $z\in \mathbb{D}$ and $V(a,b)=\left( 
\begin{array}{cc}
\bar{a} & 0 \\ 
0 & b
\end{array}
\right) $ for some $a,b\in S^{1}$. Let 
\begin{eqnarray*}
R(a,z) &=&\bar{a}L(z) \\
T(ab) &=&aV(a,b).
\end{eqnarray*}
Then for each $A\in U(1,1)$, there exists $z\in \mathbb{D}$, $a,c\in S^{1}$
such that 
\begin{equation*}
A=R(a,z)T(c)=\bar{a}\gamma _{z}\left( 
\begin{array}{cc}
1 & \bar{z} \\ 
z & 1
\end{array}
\right) \left( 
\begin{array}{cc}
1 & 0 \\ 
0 & c
\end{array}
\right) .
\end{equation*}
Set 
\begin{equation*}
S^{1}\cdot \mathcal{P}U(1,1)=\{R(a,z):a\in S^{1},\ z\in \mathbb{D}\}
\end{equation*}
and note that $T(c)\in \left\{ 1\right\} \times U(1)$. Then we have shown
that 
\begin{equation*}
U(1,1)=\left( S^{1}\cdot \mathcal{P}U(1,1)\right) (\left\{ 1\right\} \times
U(1))
\end{equation*}
is a transversal decomposition; its uniqueness follows from the uniqueness
of the polar decomposition. It is easy to show that the induced binary
operation and corresponding transversal mapping on $S^{1}\cdot \mathcal{P}%
U(1,1)$ are given by 
\begin{equation*}
R(a,z)\odot R(b,w)=R\left( ab\frac{1+\bar{z}w}{|1+\bar{z}w|},z\oplus w\right)
\end{equation*}
and 
\begin{equation*}
l(R(a,z),R(b,w))=T\left( \frac{1+\bar{z}w}{1+z\bar{w}}\right) ,
\end{equation*}
respectively, where $\oplus $ is given by (\ref{disk-op}). Straightforward
computations give 
\begin{equation*}
T(a)R(b,z)T(a)^{-1}=R(b,a^{2}z)
\end{equation*}
and 
\begin{equation*}
R(a,z)R(b,w)R(a,z)=R(ab^{2}c/\left| c\right| ,(z\oplus w)\oplus z),
\end{equation*}
where $c=1+\bar{z}w+z\bar{w}+\left| z\right| ^{2}$. Thus (G-A$_{l}$) and
(G-Bol) hold, from which it follows that $(S^{1}\cdot \mathcal{P}%
U(1,1),\odot )$ is an $A_{l}$ Bol loop. It is easy to show directly that $%
S^{1}\cdot \mathcal{P}U(1,1)$ does not satisfy AIP and hence is not a Bruck
loop.

\noindent The loop $(S^{1}\cdot \mathcal{P}U(1,1),\odot )$ is isomorphic to
a loop structure on the set $H_{\mathbb{C}}=\left\{ (x_{0},x_{1})^{t}\in 
\mathbb{C}^{2}:\left| x_{0}\right| ^{2}-\left| x_{1}\right| ^{2}=1\right\} $%
. For $(x_{0},x_{1})^{t},(y_{0},y_{1})^{t}\in H_{\mathbb{C}}$, define 
\begin{equation*}
\left( 
\begin{array}{c}
x_{0} \\ 
x_{1}
\end{array}
\right) \odot \left( 
\begin{array}{c}
y_{0} \\ 
y_{1}
\end{array}
\right) =\left( 
\begin{array}{c}
\frac{x_{0}}{\bar{x}_{0}}(\bar{x}_{0}y_{0}+\bar{x}_{1}y_{1}) \\ 
x_{0}y_{1}+x_{1}y_{0}
\end{array}
\right) .
\end{equation*}
Then $(H_{\mathbb{C}},\odot )$ is a loop. (This is a simplified version of
an example in \cite{smith-ungar}.) The following sequence of loop
homomorphisms is exact: 
\begin{equation*}
1\rightarrow S^{1}\overset{\alpha }{\rightarrow }H_{\mathbb{C}}\overset{\pi 
}{\rightarrow }\mathbb{D}\rightarrow 0,
\end{equation*}
where $\alpha (z)=(z,0)^{t}$ for $z\in S^{1}$ and $\pi \left(
(x_{0},x_{1})^{t}\right) =x_{1}/x_{0}$. In fact, $H_{\mathbb{C}}$ is a
central, invariant extension of $\mathbb{D}$ by $S^{1}$ \cite{rozga}. The
mapping $(H_{\mathbb{C}},\odot )\rightarrow (S^{1}\cdot PU(1,1),\odot
):(x_{0},x_{1})\mapsto R(x_{0}/\left| x_{0}\right| ,x_{1}/x_{0})$ is an
isomorphism of loops.
\end{example}

\begin{example}
\label{ex-proj-groups}(\emph{Projective Groups})~Let $G\leq SL(n,\mathbb{C})$
be a subgroup respecting the polar decomposition. Let $B=G\cap \mathcal{P}%
(n) $ and $H=G\cap U(n)$. Then $G=BH$ is a transversal decomposition. The
kernel of the conjugation homomorphism $U\mapsto \sigma _{U}$, where $\sigma
_{U}(A)=UAU^{\ast }$ for $A\in B$, $U\in H$, is the group $\ker (\sigma
)=G\cap \{cI:c\in \mathbb{C}\}$ of scalar matrices in $G$. Thus $PG=G/\ker
(\sigma )$ is the \emph{projective group} associated to $G$, and similarly
define $PH$. Applying (\ref{eq:factor_decomp}) to the present setting, we
have the (corefree) transversal decomposition 
\begin{equation*}
PG=B\cdot PH\text{.}
\end{equation*}
We will refer to this as a \emph{projective polar decomposition}.

As a specific example, consider the M\"{o}bius group $PSU(1,1)$. The
projective polar decomposition of this group is 
\begin{equation*}
PSU(1,1)=\mathcal{P}U(1,1)\cdot PS(U(1)\times U(1)).
\end{equation*}
The only scalar matrices in $S(U(1)\times U(1))$ are $\pm I$, and thus $%
PS(U(1)\times U(1))=S(U(1)\times U(1))/\{\pm I\}$. In terms of (\ref
{eq:H_seq}) and (\ref{eq:internal_seq}), the exact sequence of groups 
\begin{equation*}
1\rightarrow \{\pm I\}\rightarrow S(U(1)\times U(1))\rightarrow
PS(U(1)\times U(1))\rightarrow 1
\end{equation*}
induces the exact sequence of semidirect product groups 
\begin{equation*}
1\rightarrow \{\pm I\}\rightarrow SU(1,1)\rightarrow PSU(1,1)\rightarrow 1.
\end{equation*}
\end{example}

\begin{example}
\label{ex-triangular}(\emph{Upper Triangular Matrices})~Let $\mathbb{F}$ be
a field containing $1/2$. For $x,y,z\in \mathbb{F}$, let 
\begin{equation*}
T(x,y,z)=\left( 
\begin{array}{ccc}
1 & x & y \\ 
0 & 1 & z \\ 
0 & 0 & 1
\end{array}
\right) ,
\end{equation*}
and let $\mathcal{T}(3,\mathbb{F})=\{T(x,y,z):x,y,z\in \mathbb{F}\}$ be the
group of $3\times 3$ strictly upper triangular matrices over $\mathbb{F}$.
For $x_{1},x_{2}\in \mathbb{F}$, let 
\begin{equation*}
A(x_{1},x_{2})=\left( 
\begin{array}{ccc}
1 & x_{1} & \frac{1}{2}x_{1}x_{2} \\ 
0 & 1 & x_{2} \\ 
0 & 0 & 1
\end{array}
\right) ,
\end{equation*}
and let $A(3,\mathbb{F})=\{A(x_{1},x_{2}):x_{1},x_{2}\in \mathbb{F}\}$. For $%
c\in \mathbb{F}$, let 
\begin{equation*}
M(c)=\left( 
\begin{array}{ccc}
1 & 0 & c \\ 
0 & 1 & 0 \\ 
0 & 0 & 1
\end{array}
\right) ,
\end{equation*}
and let $M(3,\mathbb{F})=\{M(c):c\in \mathbb{F}\}$. Then $M(3,\mathbb{F})$
is a subgroup of $\mathcal{T}(3,\mathbb{F})$. An arbitrary matrix $%
T(x,y,z)\in \mathcal{T}(3,\mathbb{F})$ factors as follows: 
\begin{equation*}
T(x,y,z)=A(x,z)M\left( y-xz/2\right) .
\end{equation*}
It is easy to show by direct computation that this factorization is unique.
In addition, $A(x,z)=M(y-xz/2)$ if and only if $x=y=z=0$, which implies $A(3,%
\mathbb{F})\cap M(3,\mathbb{F})=\{I\}$. Thus 
\begin{equation}
\mathcal{T}(3,\mathbb{F})=A(3,\mathbb{F})\cdot M(3,\mathbb{F})
\label{triangle-decomp}
\end{equation}
is a transversal decomposition. Denote the induced binary operation (\ref
{eq:B_operation}) on $A(3,\mathbb{F})$ by $\oplus $. Then $\oplus $ is given
by 
\begin{equation}
A(x_{1},x_{2})\oplus A(y_{1},y_{2})=A(x_{1}+y_{1},x_{2}+y_{2}).
\label{A(3)-operation}
\end{equation}
It is immediate that $A(3,\mathbb{F})$ is an abelian group isomorphic to $%
\mathbb{F}^{2}$. The transversal mapping $l:A(3,\mathbb{F})\times A(3,%
\mathbb{F})\rightarrow M(3,\mathbb{F})$ is given by 
\begin{equation}
l(A(x_{1},x_{2}),A(y_{1},y_{2}))=M((x_{1}y_{2}-x_{2}y_{1})/2).
\label{triang-d}
\end{equation}
Since this is nontrivial, $A(3,\mathbb{F})$ is not a subgroup of $\mathcal{T}%
(3,\mathbb{F})$. The subgroup $M(3,\mathbb{F})$ is the center of $\mathcal{T}%
(3,\mathbb{F})$, and thus trivially normalizes $A(3,\mathbb{F})$. Thus $%
\mathcal{T}(3,\mathbb{F})$ is an internal semidirect product of the abelian
group $(A(3,\mathbb{F}),\oplus )$ with the abelian group $M(3,\mathbb{F})$.
However, as noted in Remark \ref{rem:special-cases}(4), this is \textit{not}
the usual internal semidirect product of groups.

This example turns out to generalize to $\mathcal{T}(n,\mathbb{F})$ for any $%
n\geq 3$. For $n>3$, $(A(n,\mathbb{F}),\oplus )$ is a loop but not a group,
and the homomorphism $\sigma :$ $M(n,\mathbb{F})\rightarrow \mathrm{Aut}(A(n,%
\mathbb{F}))$ is nontrivial.
\end{example}

\begin{remark}
\label{rem-sabinin-counter}Examples \ref{ex-proj-groups} and \ref
{ex-triangular} are counterexamples to a claim of Sabinin (\cite{sabinin},
Thm. 8) that (in the notation and terminology of the present paper) $\sigma
:H\rightarrow \mathrm{Sym}_{1}(B)$ is always injective.
\end{remark}

\section{External Semidirect Products}

In this section we generalize the standard semidirect product $B\rtimes H$
to the case where $H$ is not necessarily a transassociant of $\mathrm{Sym}%
_{1}(B)$, but rather there is a homomorphism $\sigma :H\rightarrow \mathrm{%
Sym}_{1}(B)$. Our discussion will show that, in a certain sense, our
definition of external semidirect product is the optimal one.

Let $B$ be a set with a distinguished element $1$ and let $H$ be a group
with identity element $e$. Assume that $B\times H$ has a binary operation
which makes it a group satisfying the following properties:

\begin{enumerate}
\item[(E1)]  $\left\{ 1\right\} \times H$ is a subgroup isomorphic to $H$.

\item[(E2)]  $(x,h)=(x,1)(1,h)$ for all $x\in B$, $h\in H$.
\end{enumerate}

\noindent Then $B\times H=(B\times \left\{ 1\right\} )\left( \left\{
1\right\} \times H\right) $ is a transversal decomposition. Indeed, by (E2), 
$(x,e)(\left\{ 1\right\} \times H)=\left\{ (x,h):h\in H\right\} $ for all $%
x\in B$, which implies that $B\times \{e\}$ is a transversal. Making the
usual identifications $B\cong B\times \left\{ 1\right\} $ and $H\cong
\left\{ 1\right\} \times H$, we have an induced operation $\cdot $ on $B$
and an induced transversal mapping $l:B\times B\rightarrow H$, both defined
by 
\begin{equation*}
(x,1)(y,1)=(x\cdot y,l(x,y))
\end{equation*}
for $x,y\in B$. We also have a mapping $m:B\times H\rightarrow H$ and (what
turns out to be) a homomorphism $\sigma :H\rightarrow \mathrm{Sym}_{1}(B)$,
both defined by 
\begin{equation*}
(1,h)(x,e)(1,h^{-1})=(\sigma _{h}(x),m(x,h))
\end{equation*}
for $x\in B$, $h\in H$. Thus given arbitrary elements $(x,h),(y,k)\in
B\times H$, we compute 
\begin{eqnarray*}
(x,h)(y,k) &=&(x,1)(1,h)(y,1)(1,h^{-1})(1,h)(1,k) \\
&=&(x,1)(\sigma _{h}(y),1)(1,m(y,h))(1,hk) \\
&=&(x\cdot \sigma _{h}(y),l(x,\sigma _{h}(y)))(1,m(y,h)hk) \\
&=&(x\cdot \sigma _{h}(y),l(x,\sigma _{h}(y))m(y,h)hk).
\end{eqnarray*}

We will use this to motivate our definition. Starting over, let $B$ be a
left loop with identity elment $1$, and let $H$ be a group with identity
element $e$. Assume there exist a mapping $l:B\times B\rightarrow H$, a
mapping $m:B\times H\rightarrow H$, and a homomorphism $\sigma :H\rightarrow 
\mathrm{Sym}_{1}(B)$. With Theorem \ref{thm:sigma_properties} as motivation,
assume the following conditions hold:

\begin{enumerate}
\item[(S1)]  For all $x,y\in B$, 
\begin{equation*}
\sigma _{l(x,y)}=L(x,y).
\end{equation*}

\item[(S2)]  For all $x\in B$, $h\in H$, 
\begin{equation*}
\sigma _{m(x,h)}=\mu _{x}(h).
\end{equation*}
\end{enumerate}

Define a binary operation on $B\times H$ by 
\begin{equation}
(x,h)(y,k)=(x\cdot \sigma _{h}(y),l(x,\sigma _{h}(y))m(y,h)hk).
\label{ext-prod}
\end{equation}
The question is thus: what are the minimal additional assumptions necessary
for $B\times H$ with the product given by (\ref{ext-prod}) to be a group?

If $B\times H$ is a group and both (E1) and (E2) hold, then $B\times
H=(B\times \left\{ e\right\} )(\left\{ 1\right\} \times H)$ is a transversal
decomposition. Thus there is an induced product on $B$, which we will denote
by $\hat{\cdot}$, an induced transversal mapping $\hat{l}:B\times
B\rightarrow H$, an induced mapping $\hat{m}:B\times H\rightarrow H$, and an
induced homomorphism $\hat{\sigma}:H\rightarrow \mathrm{Sym}_{1}(B)$. Using (%
\ref{ext-prod}), (E1) and (E2), we compute 
\begin{eqnarray*}
(x\hat{\cdot}y,\hat{l}(x,y)) &=&(x,e)(y,e) \\
&=&(x\cdot y,l(x,y)m(y,e))
\end{eqnarray*}
and 
\begin{eqnarray*}
(\hat{\sigma}_{h}(x),\hat{m}(x,h)) &=&(1,h)(x,e)(1,h^{-1}) \\
&=&(1,h)(x,h^{-1}) \\
&=&(\sigma _{h}(x),l(1,\sigma _{h}(x))m(x,h)).
\end{eqnarray*}
Thus we see that $\hat{\cdot}=\cdot $ and $\hat{\sigma}=\sigma $. In
addition, we have 
\begin{eqnarray}
\hat{l}(x,y) &=&l(x,y)m(y,e)  \label{lhat} \\
\hat{m}(x,h) &=&l(1,\sigma _{h}(x))m(x,h)  \label{mhat}
\end{eqnarray}
for all $x,y\in B$, $h\in H$. If we assume that $\hat{l}=l$ and $\hat{m}=m$,
then we have the following necessary requirements:

\begin{enumerate}
\item[(S3)]  For all $x\in B$, 
\begin{equation*}
l(1,x)=e.
\end{equation*}

\item[(S4)]  For all $x\in B$, 
\begin{equation*}
m(x,e)=e.
\end{equation*}
\end{enumerate}

Taking $y=1$ in (\ref{lhat}), applying (\ref{eq:l(1,a)=l(a,1)=1})\ to $\hat{l%
}(x,1)$, and using (S4), we obtain

\begin{enumerate}
\item[(S5)]  For all $x\in B$, 
\begin{equation*}
l(x,e)=e.
\end{equation*}
\end{enumerate}

Taking $x=1$ in (\ref{mhat}), applying (\ref{eq:m(1,h)=1}) to $\hat{m}(1,h)$%
, and using (S3), we obtain

\begin{enumerate}
\item[(S6)]  For all $h\in H$, 
\begin{equation*}
m(1,h)=e.
\end{equation*}
\end{enumerate}

Next we consider the group axioms which must be satisfied by $B\times H$. We
have 
\begin{equation*}
(x,h)(1,e)=(x\cdot \sigma _{h}(1),l(x,\sigma _{h}(1))m(1,h)h)=(x,h)
\end{equation*}
by (S5) and (S6), and 
\begin{equation*}
(1,e)(x,h)=(1\cdot \sigma _{e}(x),l(1,\sigma _{e}(x))m(x,e)h)=(x,h)
\end{equation*}
by (S3) and (S4). Therefore, the hypotheses we have so far give us that $%
(1,e)$ is the identity element of $B\times H$.

Next, we impose associativity on $B\times H$. By computing an arbitrary
product $(x,h)(y,k)(z,t)$ in two different ways, matching $H$-components
(matching $B$-components gives no new information), and simplifying, we
obtain the following technical condition which must be satisfied.

\begin{enumerate}
\item[(TC)]  For all $x,y,z\in B$, $h,k\in H$, 
\begin{multline*}
l(x\cdot \sigma _{h}(y),(L(x,\sigma _{h}(y))\mu _{y}(h)\sigma
_{hk})(z))m(z,l(x,\sigma _{h}(y))m(y,h)hk)l(x,\sigma _{h}(y))m(y,h) \\
=l(x,\sigma _{h}(y)\cdot (\mu _{y}(h)\sigma _{hk})(z))m(y\cdot \sigma
_{k}(z),h)hl(y,\sigma _{k}(z))m(z,k)h^{-1}.
\end{multline*}
\end{enumerate}

\noindent Fortunately, (TC) can be replaced by three simpler conditions to
which it is equivalent. First, taking $h=k=e$ in (TC) and using (S4), we
obtain

\begin{enumerate}
\item[(S7)]  For all $x,y,z\in B$, 
\begin{equation*}
l(x\cdot y,L(x,y)z)m(z,l(x,y))l(x,y)=l(x,y\cdot z)l(y,z).
\end{equation*}
\end{enumerate}

Second, taking $x=y=1$ in (TC) and using (S3) and (S6) (and writing $x$ for $%
z$) we obtain

\begin{enumerate}
\item[(S8)]  For all $x\in B$, $h,k\in H$, 
\begin{equation*}
m(x,hk)=m(\sigma _{k}(x),h)hm(x,k)h^{-1}.
\end{equation*}
\end{enumerate}

Finally, taking $x=1$, $k=e$ in (TC) and using (S3) and (S4) (and making the
replacements $y\rightarrow x$, $z\rightarrow y$), we obtain

\begin{enumerate}
\item[(S9)]  For all $x,y\in B$, $h\in H$, 
\begin{equation*}
l(\sigma _{h}(x),(\mu _{x}(h)\sigma _{h})(y))m(y,m(x,h)h)m(x,h)=m(x\cdot
y,h)hl(x,y)h^{-1}.
\end{equation*}
\end{enumerate}

Thus we have shown one direction of the following.

\begin{lemma}
\label{lem-tech}Condition (TC) is equivalent to conditions (S7), (S8) and
(S9).
\end{lemma}

We omit the tedious proof that (S7), (S8) and (S9) imply (TC), except to say
that starting with the left hand side of (TC), one can obtain the right hand
side by two applications of (S8), then one application of (S7), and then one
application of (S9).

Next we consider inverses in $B\times H$. Consider first an element of the
form $(x,e)$. If $(y,h)$ is to be the right inverse of $(a,e)$, then 
\begin{eqnarray*}
(1,e) &=&(x,e)(y,h) \\
&=&(x\cdot y,l(x,y)m(y,e)h).
\end{eqnarray*}
This implies $y=x^{\rho }$, and using (S4), $h=l(x,x^{\rho })^{-1}$, so that 
\begin{equation}
(x,e)^{-1}=(x^{\rho },l(x,x^{\rho })^{-1}).  \label{(a,e)-inv}
\end{equation}

Now we consider a general element $(x,h)\in B\times H$ and use (E2), (E1), (%
\ref{(a,e)-inv}) and (S3) to derive the inverse: 
\begin{eqnarray}
(x,h)^{-1} &=&((x,e)(1,h))^{-1}  \notag \\
&=&(1,h^{-1})(x^{\rho },l(x,x^{\rho })^{-1})  \notag \\
&=&(\sigma _{h^{-1}}(x^{\rho }),l(1,\sigma _{h^{-1}}(x^{\rho }))m(x^{\rho
},h^{-1})h^{-1}l(x,x^{\rho })^{-1})  \notag \\
&=&(\sigma _{h^{-1}}(x^{\rho }),m(x^{\rho },h^{-1})h^{-1}l(x,x^{\rho
})^{-1}).  \label{inv-cand}
\end{eqnarray}
We check that this candidate is indeed a right inverse: 
\begin{eqnarray*}
&&(x,h)(\sigma _{h^{-1}}(x^{\rho }),m(x^{\rho },h^{-1})h^{-1}l(x,x^{\rho
})^{-1}) \\
&=&(x\cdot x^{\rho },l(x,x^{\rho })m(\sigma _{h^{-1}}(x^{\rho
}),h)hm(x^{\rho },h^{-1})h^{-1}l(x,x^{\rho })^{-1})=(1,e).
\end{eqnarray*}
By (S8), the $H$-component of the last step of this calculation indeed
simplifies to $1$. Similar computations show that (\ref{inv-cand}) is a left
inverse, although it is not necessary to check this; a two-sided identity,
right inverses and associativity are sufficient for $B\times H$ to be a
group.

\begin{definition}
\label{defn-external}Let $B$ be a left loop and let $H$ be a group. Assume
there exist a mapping $l:B\times B\rightarrow H$, a mapping $m:B\times
H\rightarrow H$, and a homomorphism $\sigma :H\rightarrow \mathrm{Sym}%
_{1}(B) $ such that conditions (S1) through (S9) are satisfied. Define a
binary operation $\cdot $ on the set $B\times H$ by 
\begin{equation*}
(x,h)\cdot (y,k)=(x\cdot \sigma _{h}(y),l(x,\sigma _{h}(y))m(y,h)hk)
\end{equation*}
for $x,y\in B$, $h,k\in H$. Then $(B\times H,\cdot )$ is called the \emph{%
external semidirect product} of $B$ with $H$ given by $(\sigma ,l,m)$, and
is denoted $B\rtimes _{(\sigma ,l,m)}H$.
\end{definition}

Of course, the whole discussion leading up to this result was a sketch of
the proof of the following.

\begin{theorem}
\label{thm-external-group}$(B\rtimes _{(\sigma ,l,m)}H,\cdot )$ is a group.
\end{theorem}

\begin{remark}
\label{rem-ext}As with the internal product, there are various special cases
of the external product which are of interest.

\begin{enumerate}
\item  If the transversal mapping $l:B\times B\rightarrow H$ can be chosen
to be trivial, i.e., $l(x,y)=e$ for all $x,y\in B$, then $B$ is a group (by
(S1)). In this case, we have an external version of Jajcay's ``rotary
product'' of groups \cite{jajcay}.

\item  If the mapping $m:B\times B\rightarrow H$ can be chosen to be
trivial, i.e., $m(x,h)=e$ for all $x\in B$, $h\in H$, then $B$ is an $A_{l}$%
-loop (by (S2)).

\item  If both $l$ and $m$ can be chosen to be trivial, then $B\rtimes
_{(\sigma ,l,m)}H=B\rtimes _{\sigma }H$ is the usual external semidirect
product of groups.

\item  If $\sigma (H)=\{I\}$, then $B$ is a group, but if $l:B\times
B\rightarrow H$ is nontrivial, then $B$ is not (isomorphic to) a subgroup of 
$B\rtimes _{(\sigma ,l,m)}H$.
\end{enumerate}
\end{remark}

Our examples are $A_{l}$ Bol loops. Thus we will choose $m$ to be trivial,
so that (S2), (S4), (S6) and (S8) are trivial. (S7) simplifies to

\begin{enumerate}
\item[(S7')]  For all $x,y,z\in B$, $h\in H$, 
\begin{equation*}
l(x\cdot y,L(x,y)z)l(x,y)=l(x,y\cdot z)l(y,z).
\end{equation*}
\end{enumerate}

and (S9) simplifies to

\begin{enumerate}
\item[(S9')]  For all $x,y\in B$, $h\in H$, 
\begin{equation*}
l(\sigma _{h}(x),\sigma _{h}(y))=hl(x,y)h^{-1}.
\end{equation*}
\end{enumerate}

\begin{example}
The mapping $\phi :S^{1}\rightarrow S^{1}$ defined by $\phi (a)=a^{2}$ is a
homomorphism of groups, and we consider the target copy of $S^{1}$ to be a
subgroup of $\mathrm{Aut}(\mathbb{D},\oplus )$. Define $l:\mathbb{D}\times 
\mathbb{D}\rightarrow S^{1}$ by 
\begin{equation*}
l(x,y)=\frac{1+x\bar{y}}{|1+\bar{x}y|}
\end{equation*}
for $x,y\in \mathbb{D}$. Then $\phi (l(x,y))=L(x,y)=(1+x\bar{y})/(1+\bar{x}%
y) $, as in Example \ref{ex-disk}. Define $m:\mathbb{D}\times
S^{1}\rightarrow S^{1}$ trivially. It is easy to check that (S1), (S3),
(S6), (S7') and (S9') are satisfied. Thus we have an external semidirect
product $\mathbb{D}\rtimes _{(\phi ,l)}S^{1}$. The exact sequence of groups 
\begin{equation*}
1\rightarrow \{\pm 1\}\rightarrow S^{1}\overset{\phi }{\rightarrow }%
S^{1}\rightarrow 1
\end{equation*}
induces an exact sequence of semidirect product groups 
\begin{equation*}
1\rightarrow \{\pm 1\}\rightarrow \mathbb{D}\rtimes _{(\phi ,l)}S^{1}%
\overset{\hat{\phi}}{\rightarrow }\mathbb{D}\rtimes S^{1}\rightarrow 1
\end{equation*}
where $\hat{\phi}(z,a)=(z,\phi (a))=(z,a^{2})$. \smallskip We now show the
relationship between this construction and Example \ref{ex-proj-groups}.
Recall the isomorphism $Q:S^{1}\rightarrow S(U(1)\times U(1))$ given by $%
Q(a)=\left( 
\begin{array}{cc}
\bar{a} & 0 \\ 
0 & a
\end{array}
\right) $. As the following diagram indicates, there exists an isomorphism $%
\hat{Q}:S^{1}\rightarrow PS(U(1)\times U(1))$. 
\begin{equation*}
\begin{array}{ccccccccc}
1 & \rightarrow & \{\pm 1\} & \rightarrow & S^{1} & \overset{\phi }{%
\rightarrow } & S^{1} & \rightarrow & 1 \\ 
&  & \downarrow &  & \downarrow _{Q} &  & \downarrow _{\hat{Q}} &  &  \\ 
1 & \rightarrow & \{\pm I\} & \rightarrow & S(U(1)\times U(1)) & \rightarrow
& PS(U(1)\times U(1)) & \rightarrow & 1
\end{array}
\end{equation*}
This is explicitly given by 
\begin{equation*}
\hat{Q}(a)=\left\{ 
\begin{array}{cc}
\left[ Q\left( \frac{1+a}{|1+a|}\right) \right] & \text{if }a\neq -1 \\ 
\lbrack Q(i)] & \text{if }a=-1
\end{array}
\right.
\end{equation*}
where the square brackets denote equivalence classes. (Note that $%
((1+a)/|1+a|)^{2}=a$.) On the other hand, the mapping $(z,a)\mapsto P(z)Q(a)$
gives an isomorphism from $\mathbb{D}\rtimes _{\phi }S^{1}$ to $SU(1,1)=%
\mathcal{P}U(1,1)\cdot S(U(1)\times U(1))$. This and the following diagram
yield an isomorphism between $\mathbb{D}\rtimes S^{1}$ and $PSU(1,1)$: 
\begin{equation*}
\begin{array}{ccccccccc}
1 & \rightarrow & \{(0,\pm 1)\} & \rightarrow & \mathbb{D}\rtimes _{(\phi
,l)}S^{1} & \overset{\hat{\phi}}{\rightarrow } & \mathbb{D}\rtimes S^{1} & 
\rightarrow & 1 \\ 
&  & \downarrow &  & \downarrow _{PQ} &  & \downarrow _{P\hat{Q}} &  &  \\ 
1 & \rightarrow & \{\pm I\} & \rightarrow & SU(1,1) & \rightarrow & PSU(1,1)
& \rightarrow & 1
\end{array}
\end{equation*}
The isomorphism is explicitly given by $(z,a)\mapsto P(z)\hat{Q}(a)$.
\end{example}

\begin{example}
Let $\mathbb{F}$ be a field containing $1/2$. Let $P=\mathbb{F}^{2}$ with
the operation of vector addition, and let $H=\mathbb{F}$ with the operation
of addition. Define $l:\mathbb{F}^{2}\times \mathbb{F}^{2}\rightarrow 
\mathbb{F}$ by $l((x_{1},x_{2}),(y_{1},y_{2}))=\frac{1}{2}%
(x_{1}y_{2}-x_{2}y_{1})$. Define $\phi :\mathbb{F}\rightarrow \mathrm{Aut}(%
\mathbb{F}^{2})$ trivially: $\phi (c)=I$. Finally, define $m:\mathbb{F}%
^{2}\times \mathrm{Aut}(\mathbb{F}^{2})\rightarrow \mathrm{Aut}(\mathbb{F}%
^{2})$ trivially. Then (S1), (S3), (S6), (S7') and (S9') are satisfied. We
thus have an external semidirect product $\mathbb{F}^{3}=\mathbb{F}%
^{2}\rtimes _{(\phi ,l)}\mathbb{F}$. The operation is given by 
\begin{equation*}
(x_{1},x_{2},x_{3})\cdot
(y_{1},y_{2},y_{3})=(x_{1}+y_{1},x_{2}+y_{2},x_{3}+y_{3}+\frac{1}{2}%
(x_{1}y_{2}-x_{2}y_{1})).
\end{equation*}
As in Remark \ref{rem-ext}(4), this external semidirect product of groups
does not reduce to the usual one. Clearly this external semidirect product
is isomorphic to the internal semidirect product of Example \ref
{ex-triangular}.
\end{example}

We conclude with a few remarks about the relationship of the external
semidirect product with the standard and internal semidirect products. These
naturally generalize the usual relationships between these products of
groups. First, if $G=BH$ is an internal semidirect product of a left loop $%
(B,\cdot )$ with the subgroup $H$, then the mapping $g=xh\mapsto (x,h)$ is
clearly an isomorphism of $BH$ with $B\rtimes _{(\sigma ,l,m)}H$. On the
other hand, by factoring a given external semidirect product as $B\rtimes
_{(\sigma ,d,m)}H\cong (B\times \{e\})(\{1\}\times H)$, it is easy to see
that it is isomorphic to an internal semidirect product; this was
essentially our starting point for deriving the definition of external
semidirect product. The standard semidirect product is, of course, a special
case of the external semidirect product. On the other hand, if $B\rtimes
_{(\sigma ,l,m)}H$ is an external semidirect product, then the natural
mapping $\hat{\sigma}:B\rtimes _{(\sigma ,d,m)}H\rightarrow B\rtimes \sigma
(H):(x,h)\mapsto (x,\sigma _{h})$ is an epimorphism. We have $\ker (\hat{%
\sigma})=\{1\}\times \ker (\sigma )$, and thus the exact sequence of groups 
\begin{equation*}
1\rightarrow \ker (\sigma )\rightarrow H\rightarrow \sigma (H)\rightarrow 1
\end{equation*}
induces an exact sequence of semidirect product groups 
\begin{equation*}
1\rightarrow \ker (\hat{\sigma})\rightarrow B\rtimes _{(\sigma
,l,m)}H\rightarrow B\rtimes \sigma (H)\rightarrow 1\text{.}
\end{equation*}

\end{document}